\documentclass{amsart}

\usepackage{amssymb}
\usepackage[all]{xy}

\newtheorem{thm}{Theorem}

\newtheorem{lem}[thm]{Lemma}
\newtheorem{cor}[thm]{Corollary}

\newtheorem{prop}[thm]{Proposition}

\theoremstyle{definition}
\newtheorem{defn}[thm]{Definition}

\newtheorem{say}[thm]{}
\newtheorem{exmp}[thm]{Example}

\newtheorem{ques}[thm]{Question}    

\newtheorem*{ack}{Acknowledgments}      

\newtheorem{defn-thm}[thm]{Definition--Theorem}  
\newtheorem{defn-lem}[thm]{Definition--Lemma}  

\theoremstyle{remark}

\renewcommand{\o}[0]{{\mathcal O}} 
\newcommand{\z}[0]{{\mathbb Z}}

\renewcommand{\r}[0]{{\mathbb R}} 

\renewcommand{\a}[0]{{\mathbb A}}

\newcommand{\p}[0]{{\mathbb P}}

\newcommand{\q}[0]{{\mathbb Q}}
\newcommand{\map}[0]{\dasharrow}
\newcommand{\qtq}[1]{\quad\mbox{#1}\quad}

\newcommand{\supp}[0]{\operatorname{Supp}}    
    
\newcommand{\codim}[0]{\operatorname{codim}}    
    
\newcommand{\proj}[0]{\operatorname{Proj}}

\newcommand{\ex}[0]{\operatorname{Ex}}    
\newcommand{\diff}[0]{\operatorname{Diff}}

\newcommand{\cl}[0]{\operatorname{Cl}}

\newcommand{\rup}[1]{\lceil{#1}\rceil}
\newcommand{\rdown}[1]{\lfloor{#1}\rfloor}

\newcommand{\simq}[0]{\sim_{\q}}

\newcommand{\simr}[0]{\sim_{\r}}

\newcommand{\depth}[0]{\operatorname{depth}} 
\newcommand{\tsum}[0]{\textstyle{\sum}}

\newcommand{\coeff}[0]{\operatorname{coeff}}

\def\into{\DOTSB\lhook\joinrel\to}

\def\loccoh#1.#2.#3.#4.{H^{#1}_{#2}(#3,#4)}

\DeclareMathAlphabet{\mathchanc}{OT1}{pzc}%
                                {m}{it}

\usepackage[all]{xy}\xyoption{dvips}

\begin{document}
\bibliographystyle{amsalpha}

 \title{Log-plurigenera in stable families}
 \author{J\'anos Koll\'ar}

\begin{abstract} We study the flatness of log-pluricanonical sheaves on stable families of varieties.
\end{abstract}

 \maketitle

A key insight of \cite{ksb} is that,  in dimension at least 2, the correct objects  of moduli theory   are   flat morphisms $f:X\to S$ 
whose fibers  are  varieties with log canonical singularities and 
such that  $mK_{X/S}$ is Cartier for some $m=m(X, S)>0$.  The  latter assumption is not always easy to understand, but, if $S$ is reduced then  it is equivalent to the condition
\begin{itemize}
\item   $\omega_{X/S}^{[m]}$ is flat over $S$ and commutes with base change for every $m\in \z$. 
\end{itemize}

In studying the moduli of  pairs, the right concept is less clear. One should consider morphisms $f:(X,\Delta)\to S$ such that $f:X\to S$ is flat, all fibers $(X_s, \Delta_s)$ are    semi-log-canonical and $K_{X/S}+\Delta$ is $\r$-Cartier. Such morphisms are called {\it locally stable;}
see \cite{k-modsurv} for a survey and \cite{k-modbook} for   a detailed treatment.

However, these conditions are not sufficient; there are problems especially over non-reduced bases \cite{k-alt} and in positive characteristic \cite[14.7]{hac-kov}. One difficulty is that the sheaf variant of the Cartier assumption is not well understood.

The best log-analog of  $\omega_{X/S}^{[m]}$ is the twisted version $\omega_{X/S}^{[m]}\bigl(\rdown{m\Delta}\bigr) $. 
 The main theorem of this note shows that these sheaves also behave well, provided   the coefficients of  $\Delta$ are not too small. 

For a divisor $\Delta=\sum_{i\in I} a_i D_i$, where the $D_i$ are distinct prime divisors, we write $\coeff\Delta:=\{a_i :i\in I\}$. In the semi-log-canonical cases $\coeff\Delta\subset [0, 1]$.  
We set  $\rdown{\Delta}:=\sum_{i\in I} \rdown{a_i} D_i$,
$\rup{\Delta}:=\sum_{i\in I} \rup{a_i} D_i$ and
$\{\Delta\}:=\sum_{i\in I} \{a_i\} D_i=\Delta-\rdown{\Delta}$.

\begin{thm}  \label{half.pgsh.inv.thm}
Let $S$ be a reduced scheme over a field of characteristic 0 and $f:(X,\Delta) \to S$ a locally stable morphism with normal generic fibers. Assume  that $\coeff\Delta\subset [\frac12, 1]$.
Then, for every $m\in\z$, the sheaves
$$
\omega_{X/S}^{[m]}\bigl(\rdown{m\Delta}\bigr)
$$
are flat over $S$ and commute with base change (cf.\ Definition~\ref{comm.base.ch.defn}).
\end{thm}

{\it Warning \ref{half.pgsh.inv.thm}.1.} If some of the coefficients equal $\frac12$, we have to be  careful about what we mean by $\rdown{m\Delta}$ and  commuting with base change; see Paragraph~\ref{half.rest.well.lem} for a general discussion and Paragraph~\ref{1/2.prob.say} for the coefficient  $\frac12$ case.
\medskip

The assumption on having normal generic fibers is probably superfluous, see Question \ref{what.about.slc.q}. Aside from this,
the theorem seems quite sharp.  In  (\ref{half.main.exmp}.6) we give  examples  
such that $\coeff\Delta$ is arbitrarily close to  $\frac12$ yet
$\omega_{X/S}^{[2]}\bigl(\rdown{2\Delta}\bigr)$ does not commute with base change.
Also, for $m\neq r$ the sheaves
$\omega_{X/S}^{[m]}\bigl(\rdown{r\Delta}\bigr)$ 
usually do not have similar properties, see  (\ref{half.main.exmp}.5).
In example (\ref{half.main.exmp}.4)  the pluricanonical sheaf
 $ \omega_{X/S}^{[m]}$ commutes with base change only for $m=0$ and  $m=1$.

In Theorem~\ref{half.pgsh.inv.thm} we allow $\Delta$ to be an $\r$-divisor.
In this case $\omega_{X/S}^{[m]}\bigl(\rdown{m\Delta}\bigr)$ need not be locally free for any $m\neq 0$; see (\ref{half.main.exmp}.8). 

In Section \ref{sec.0} we state some corollaries and variants while in  Section \ref{sec.0.b} we recall some relevant facts of the moduli theory if pairs. 
In Section \ref{sec.1} we study small modifications of lc pairs in general.  These show us how to transform divisors like $mK_X+\rdown{m\Delta} $ into a  $\q$-Cartier divisor in an economic way. Usually this can be done in several ways and we pin down some especially good choices. Most of these results do not hold for slc pairs by \cite[1.40]{kk-singbook}; see also  Example~\ref{1.40.kk-singbook.exmp}.

In Section \ref{sec.2} we prove Proposition~\ref{tech.prop.2} which is a generalization of Proposition~\ref{tech.prop.0} and then  
in Section \ref{sec.3} we establish Corollary \ref{st.pg.inv.cor}. We end by a collection of relevant examples in Section \ref{sec.4}.  

\begin{ack} I thank   A.J.\ de~Jong and Chenyang~Xu for insightful comments. 
Partial  financial support    was provided  by  the NSF under grant number
 DMS-1362960.
\end{ack}

\section{Consequences and variants}\label{sec.0}

The first corollary of  Theorem~\ref{half.pgsh.inv.thm} is the
deformation invariance of the Hilbert function  for certain locally stable morphisms. (Note that Warning~\ref{half.pgsh.inv.thm}.1 also applies to the 
 results in this section;  see Paragraph~\ref{1/2.prob.say} if some of the coefficients of the boundary divisor  equal $\frac12$.)

\begin{cor} \label{st.Hf.inv.cor}
 Let $S$ be a  connected  scheme over a field of characteristic 0 and $f:(X,\Delta) \to S$ a proper, locally stable morphism  with normal generic fibers such that $\coeff\Delta\subset [\frac12, 1]$. Then the Hilbert function of the fibers
$$
\chi\bigl(X_s, \omega_{X_s}^{[m]}\bigl(\rdown{m\Delta_s}\bigr)\bigr)
\qtq{is  independent of  $s\in S$.} \qed
$$
\end{cor}

In general the Hilbert function considered above  is not a polynomial in $m$, but,
if $r(K_{X/S}+\Delta)$ is Cartier, then  it can be written as a polynomial whose coefficients are periodic functions in $m$ with period  $r$. 
If we assume in addition that $K_{X/S}+\Delta$ is $f$-ample then, by Serre vanishing, 
$$
\chi\bigl(X_s, \omega_{X_s}^{[m]}\bigl(\rdown{m\Delta_s}\bigr)\bigr)=
H^0\bigl(X_s, \omega_{X_s}^{[m]}\bigl(\rdown{m\Delta_s}\bigr)\bigr)
\qtq{for} m\gg 1,
$$
but it is not clear for which values of $m$ does this hold. 
However we get the optimal 
deformation invariance of plurigenera if we restrict the coefficients further.

\begin{cor} \label{st.pg.inv.cor}
 Let $S$ be a reduced scheme over a field of characteristic 0 and $f:(X,\Delta) \to S$ a stable morphism  with normal generic fibers such that $\coeff\Delta\subset \{\frac12, \frac23, \frac34, \dots, 1\}$. Then, for every $m\geq 2$,
\begin{enumerate}
\item  $R^if_*\omega_{X/S}^{[m]}\bigl(\rdown{m\Delta}\bigr)=0$ for $i>0$ and
\item $f_*\omega_{X/S}^{[m]}\bigl(\rdown{m\Delta}\bigr)$ is locally free and commutes with base change.
\end{enumerate}
\end{cor}

Examples~\ref{c.near.1.h.coh.exmp}--\ref{c.near.1.h.coh.exmp.4} show that (\ref{st.pg.inv.cor}.1) can fail for  any other coefficient. 
It is, however, possible that if we fix the relative dimension and
$\coeff\Delta$ is close enough to 1 then (\ref{st.pg.inv.cor}.1) holds.
It would be good to know whether (\ref{st.pg.inv.cor}.2) holds more generally.

Rather general arguments involving hulls and husks reduce the proof of Theorem~\ref{half.pgsh.inv.thm} to the cases when $S$ is  smooth of dimension 1; see the proof of Theorem~\ref{half.pgsh.inv.D.thm} and Proposition~\ref{comm.base.ch.prop}. In these cases  $(X, \Delta)$ itself is an slc pair, and then Theorem~\ref{half.pgsh.inv.thm} can be reformulated as a version of Serre's $S_3$ property.

\begin{defn} \label{S3.on.S.defn} Recall that a sheaf $F$ on a  scheme $X$ is called $S_m$ if $$
\depth_xF\geq \min\{m, \dim_xF\}\quad \forall\ x\in X.
$$  As a slight variant, we say
that $F$ is {\it  $S_m$ along a subscheme} $Z\subset X$ if $\depth_xF\geq \min\{m, \dim_xF\}$ holds for every   $x\in Z$.

We are mostly interested in the $S_3$ condition. Note that if
$H\subset X$ is a Cartier divisor whose defining equation is  not a zero-divisor on $F$ then $F$ is $S_3$ along $H$ iff  $F|_H$ is $S_2$.
\end{defn}

\begin{prop} \label{tech.prop.0}
Let $(X,H+\Delta)$ be an lc pair  over a field of characteristic 0 where $H$ is Cartier and $\coeff\Delta\subset [\frac12, 1]$. 
Then $\omega_X^{[m]}\bigl(\rdown{m\Delta}-D\bigr)$ is $S_3$ along $H$
for every $m\in \z$ and every divisor $D\subset \rdown{\Delta}$.
\end{prop}

If, by happenstance, $mK_X+\rdown{m\Delta} $ is $\q$-Cartier, then 
this   follows from  \cite{ale-lim, k-dep};  
see also  \cite[7.20]{kk-singbook} and (\ref{alex.kol.depth.3}). Thus our main focus is on the cases when
 $mK_X+\rdown{m\Delta} $ is not $\q$-Cartier.
Having an extra divisor $D$ in Proposition~\ref{tech.prop.0} allows us to prove a stronger version of  Theorem   \ref{half.pgsh.inv.thm}.

\begin{thm}  \label{half.pgsh.inv.D.thm}
Let $S$ be a reduced scheme over a field of characteristic 0 and $f:(X,\Delta) \to S$ a locally stable morphism with normal generic fibers. Assume  that $\coeff\Delta\subset [\frac12, 1]$ and let  $D\subset \rdown{\Delta}$.
Then, for every $m\in\z$, the sheaves
$$
\omega_{X/S}^{[m]}\bigl(\rdown{m\Delta}-D\bigr)
$$
are flat over $S$ and commute with base change.
\end{thm}

Proof of Proposition~\ref{tech.prop.0} $\Rightarrow$  Theorem~\ref{half.pgsh.inv.D.thm} $\Rightarrow$ Theorem~\ref{half.pgsh.inv.thm}.
Setting $D=0$ shows that Theorem~\ref{half.pgsh.inv.D.thm} $\Rightarrow$ Theorem~\ref{half.pgsh.inv.thm}. In order to prove the first implication, we
 use the theory of {\it hulls and husks;} see
\cite{k-hh} or \cite[Chap.9]{k-modbook}.


Fix $m\in \z$. As we note in   \cite[4.26]{k-modbook}, 
there is a closed subset $Z\subset X$ such that $K_{X/S}$ and 
$\rdown{m\Delta}-D$ are Cartier on $X\setminus Z$ and 
$Z\cap X_s$ has codimension $\geq 2$ in $X_s$ for every $s\in S$. 
 
We aim to apply  Proposition~\ref{comm.base.ch.prop}  to
$U:=X\setminus Z$ with injection $j:U\into X$ and 
$F:=\omega_{U/S}^{[m]}(\rdown{m\Delta|_U}-D|_U)$. 

Lemma~\ref{comm.base.ch.lem} says that Proposition~\ref{tech.prop.0}  is equivalent to the assumption (\ref{comm.base.ch.prop}.2). Thus we get
that (\ref{comm.base.ch.prop}.1) also holds, and the latter  is just a reformulation of the claim of Theorem~\ref{half.pgsh.inv.D.thm}. \qed



\medskip

There are several results \cite{ale-lim, k-dep} that  guarantee that certain divisorial sheaves are  Cohen--Macaulay  (abbreviated as CM)  or at least $S_3$. We recall these in 
Theorem~\ref{alex.kol.depth.3}, see also 
\cite[7.20]{kk-singbook} and \cite[Sec.7.1]{fuj-book} for detailed treatments.
The following variant of Proposition~\ref{tech.prop.0} is closely related to them.

\begin{prop} \label{tech.prop.B}
Let $(X,H+\Delta)$ be an lc pair  over a field of characteristic 0 where $H$ is $\q$-Cartier and $\coeff\Delta\subset (\frac12, 1]$.    Let $B$ be a Weil $\z$-divisor  and $\Delta_3$ an effective $\r$-divisor such that 
 $B\simr -\Delta_3$,  $\Delta_3\leq \rup{\Delta}$ and  $\rdown{\Delta_3}\leq \rdown{\Delta}$.

Then $\o_X(B)$ is $S_3$ along $H$.
\end{prop}

A stronger version, allowing some coefficients to be $\frac12$, is proved in 
(\ref{tech.prop.2}).

The main open question is the following.

\begin{ques}\label{what.about.slc.q} What happens for non-normal slc pairs? 

Note that we derive Propositions~\ref{tech.prop.0} and \ref{tech.prop.B}---and hence also Theorems~\ref{half.pgsh.inv.thm} and \ref{half.pgsh.inv.D.thm}---using Proposition~\ref{pot.lc.can.mod}, which asserts that  certain small modifications of $X$ exist. The analogous small modifications need not exists for slc pairs; see
\cite[1.40]{kk-singbook} or  Example~\ref{1.40.kk-singbook.exmp}. 
However,   Proposition~\ref{tech.prop.B} and  Theorem~\ref{half.pgsh.inv.thm} frequently hold even when Proposition~\ref{pot.lc.can.mod} fails.

Some non-normal cases are treated in \cite{k-lpg2}. 
\end{ques}

The divisor $H$ does not play any role in the conclusion of  Proposition~\ref{tech.prop.0}, but it restricts the possible choices of $(X, \Delta)$. Example~\ref{2/3.v2.exmp} shows that $H$ is necessary, but I have no counter example to the following variant.  

\begin{ques} \label{tech.prop.B.ques}
Let $(X,\Delta)$ be an slc pair  over a field of characteristic 0 where  $\coeff\Delta\subset [\frac23, 1]$.   Let  $x\in X$ be a point of codimension $\geq 3$ that is not an lc center of $(X,\Delta)$. 
Is it true that  $\depth_x\omega_{X}^{[m]}\bigl(\rdown{m\Delta}\bigr)\geq 3$?
\end{ques}

While the Serre dual of a CM   sheaf is  CM (cf.\ \cite[5.70]{km-book}), the  Serre dual of an $S_3$ sheaf need not be  $S_3$; see \cite[1.6]{patakfalvi}. Thus it is not clear whether
the dual versions of our results also hold. 
The Serre dual of    $\omega_X^{[m]}\bigl(\rdown{m\Delta}-D\bigr)$ is  $\omega_X^{[1-m]}\bigl(-\rdown{m\Delta}+D\bigr)$. Changing $m$ to $-m$
we get the following. 

\begin{ques}\label{SD.sheaf.ques.2}  Using the notation and assumptions of
Proposition~\ref{tech.prop.0}, is the sheaf $$\omega_X^{[m+1]}\bigl(\rup{m\Delta}+D\bigr)\qtq{$S_3$ along $H$?}
$$
\end{ques}

\begin{say}[Method of proof]\label{S3.method.say}
The idea is similar to the ones used in
\cite{k-dep, k-1/6}.
 
 Let $g:Y\to X$ be a proper morphism of normal varieties, $F$ a coherent sheaf on $X$, $H\subset X$ a Cartier divisor and $H_Y:=g^*H$.  Assuming that $F$ is $S_m$ along $H_Y$, we   would like to understand when
$g_*F$ is $S_m$ along $H$.
If (the local equation of)  $H_Y$ is not a zero divisor on $F$  then the sequence
$$
0\to F(-H_Y)\to F\to F|_{H_Y}\to 0
\eqno{(\ref{S3.method.say}.1)}
$$
is exact. By push-forward we get the exact sequence
$$
0\to g_*F(-H)\to g_*F\to g_*\bigl(F|_{H_Y}\bigr)\to R^1g_*F(-H_Y)\cong \o_X(-H)\otimes R^1g_*F
\eqno{(\ref{S3.method.say}.2)}
$$
Thus we see that $g_*F$ is $S_m$ along $H$ if
\begin{enumerate}\setcounter{enumi}{2}
\item[(3.a)]  $R^1g_*F=0$ and 
\item[(3.b)]   $g_*\bigl(F|_{H_Y}\bigr)$ is  $S_{m-1}$ along $H$.
\end{enumerate}
In  many cases, for instance if $g$ is an isomorphism outside $H_Y$,  these conditions are also necessary.

Our main interest is in the cases when $F$ is a divisorial sheaf. 
Using a Kodaira-type vanishing theorem, 
(3.a) needs some positivity condition on $F$. By contrast,  we see in Lemma~\ref{push.wfd.lem} that  (3.b) needs some negativity  condition on $F$.

In general one can not satisfy both of these restrictions, but choosing $Y$ carefully and varying the boundary divisor gives some  wiggle room.
\end{say}

\section{Comments on the moduli of pairs}\label{sec.0.b}

The general theory of stable and locally stable maps is treated in
\cite{k-modbook}, see also \cite{k-modsurv}. Here we discuss one special aspect of it, the definition of the divisoral part of the fiber. This explains why the condition
$\coeff\Delta\subset [\frac12, 1]$ is necessary and we also need later some of its easy but potentially confusing properties.

\begin{say}[Restriction and rounding down] \label{half.rest.well.lem}
Let $f: (X, \Delta)\to S$ be a locally stable morphism. Here we consider the following problem.
\medskip

{\it Question \ref{half.rest.well.lem}.1.}
Given a point $s\in S$, 
how can we compare the divisor-fibers $\rdown{m\Delta_s}$ and $\rdown{m\Delta}_s$?
\medskip

More generally, let 
$\Theta$ be  an  $\r$-divisor such that 
  $\supp \Theta\subset \supp \Delta$. What is the relationship between  
$\rdown{\Theta_s}$ and $\rdown{\Theta}_s$?

The definition of local stability  in
\cite[Chap.3]{k-modbook} 
is set up so that  $\Delta_s$ is defined using base change  to the spectrum of a DVR  $(t, T)\to S$. In this case $X_t$ is a Cartier divsor in $X_T$, thus we consider the folowing more general variant. 
\medskip

{\it Question \ref{half.rest.well.lem}.2.}
Let $(X,H+\Delta)$ be an lc pair where $H$ is a  Cartier divisor. 
Let  $\Theta$ be  an  $\r$-divisor such that 
  $\supp \Theta\subset \supp \Delta$. What is the relationship between  
$\rdown{\Theta|_H}$ and $\rdown{\Theta}|_H$?
\medskip

In order to answer this,  let $W\subset H$ be an irreducible component of $H\cap\supp \Theta$.
Then $W$ is also an irreducible component of $H\cap\supp\Delta$.  Since 
$(X,H+\Delta)$ is lc, this implies that $X$ and $H$ are smooth at the generic point of $W$. Write $\Delta=\sum a_i D_i$ and let $m_i$ denote the intersection multiplicity of $H$ and $D_i$ along $W$.  Then $\coeff_W(\Delta|_H)=\sum_i m_ia_i$.   Since 
$(X,H+\Delta)$ is lc, so is $(H, \Delta|_H)$, hence $ \sum_i m_ia_i\leq 1$. 
This inequality does not cary much information about the small $a_i$, but it becomes stronger as the $a_i$ become larger. In particular, 
 at most one of the $a_i>\frac12$ can appear with nonzero multiplicity $m_i$ and then necessarily $m_i=1$.  This proves  the next claim.
\medskip

{\it Claim \ref{half.rest.well.lem}.3.}
Assume also that 
 $\supp \Theta\subset \supp \bigl(\Delta^{>1/2}\bigr)$. Then
\begin{enumerate}
\item[(a)]   $\coeff(\Theta|_H)\subset \coeff\Theta$ and
\item[(b)]  $\rdown{\Theta|_H}=\rdown{\Theta}|_H$. 
 \qed
\end{enumerate}
\medskip

If  $W$ is not an  lc center of  $(X,H+\Delta)$ then
we get the stronger inequality  $ \sum_i m_ia_i< 1$, and the above conclusions also apply if we allow  $a_i= \frac12$. 

\medskip

{\it Claim \ref{half.rest.well.lem}.4.} The conclusions
(\ref{half.rest.well.lem}.3.a--b) also hold if $\supp \Theta\subset \supp \bigl(\Delta^{\geq 1/2}\bigr)$ and none of the  codimension 2 lc centers of  $(X,H+\Delta)$ is contained in $H$. \qed




\end{say}

Next we discuss what happens if some of the coefficients equal $\frac12$.

\begin{say}[Problems with coefficient $\frac12$] \label{1/2.prob.say} We have to be  careful with bookkeeping if $\Delta$ contains divisors with coefficient $=\frac12$. This is best illustrated by some examples.

(\ref{1/2.prob.say}.1)  Consider  $\pi: \a^2_{xu}\to \a^1_u$ with
$\Delta=\frac12(x-u=0)+ \frac12(x+u=0)$. The generic fiber is a line with 2 points with weight $\frac12$, the special fiber over $(u=0)$ is  a line with 1 point with weight $1$. The solution is that we remember that the boundary on the special fiber is
$\Delta_0:=\frac12(x=0)+ \frac12(x=0)$ and so we declare that
$\rdown{\Delta_0}:=\rdown{\frac12}(x=0)+\rdown{\frac12}(x=0)=0$.

(\ref{1/2.prob.say}.2)  Consider  $\pi: \a^2_{xu}\to \a^1_u$ with
$\Delta=\frac12(x^2-u=0)$.  Again  the special fiber is  a line with 1 point with weight $1$.  The solution is that we write the the special fiber as
$\Delta_0:=\frac12\bigl(2(x=0)\bigr)$ and so we declare that
$\rdown{\Delta_0}:=\rdown{\frac12}\bigl(2(x=0)\bigr)=0$.
\medskip

These two examples describe what happens over a normal base scheme $S$.
In this case we write  $\Delta=\sum a_i D^i$ where the $D^i$ are prime divisors. 
If $T\to S$ is any base change then we write
$\Delta_T=\sum a_i D^i_T$. Even though the $D^i_T$ need not be prime divisors,
we set
$$
\rdown{m\Delta_T}:=\tsum \rdown{ma_i}D^i_T.
\eqno{(\ref{1/2.prob.say}.3)}
$$
If $a_i>\frac12$ then $D^i_T$ has no multiple  components and no  irreducible components in common  with any of the other $D^j_T$ by (\ref{half.rest.well.lem}.3).  Thus the only change is on how we count those divisors that have coefficient $\frac12$.

Over reducible bases, other complication can arise. 

(\ref{1/2.prob.say}.4) Consider  $X:=(uv=0)\subset \a^3_{xuv}$ and 
$S:=(uv=0)\subset \a^2_{uv}$. Set 
$\Delta:=\frac12(x-u=v=0)+ \frac12(x+u=v=0)+(x=u=0)$. 
Here we can not view $\rdown{\Delta}$ as a family of divisors in any sensible way. Along the $u$-axis we should get   $\rdown{\Delta}=0$ but along the $v$-axis $\rdown{\Delta}=(x=u=0)$. These can not be reconciled over the origin. 

(\ref{1/2.prob.say}.5)  The solution is to work with morphisms  $f:(X,\Delta=\sum a_i D_i) \to S$
where each $D_i$ is a  $\z$-divisor on $X$ that is Cartier at the generic points of $X_s\cap \supp D_i$ for every $s\in S$.  (This condition is automatic if $S$ is normal  by  \cite[4.2]{k-modbook}.)
This method is formalized by the concept of {\it marked pairs} defined in   \cite[Sec.4.6]{k-modbook}
\end{say}

Next we discuss some results that were used in the proof of 
Theorem~\ref{half.pgsh.inv.D.thm}.

\begin{defn} \label{comm.base.ch.defn}  Let   $f:X\to S$ a  morphism and   $Z\subset X$ be a closed subset such that
$Z\cap X_s$ has codimension $\geq 2$ in $X_s$ for every $s\in S$. 
Set $U:=X\setminus Z$ with injection $j:U\into X$ and let $F$ be a coherent sheaf on $U$. For every morphism  $g:T\to S$ we have a base-change diagram
$$
\begin{array}{ccccc}
U_T & \stackrel{j_T}{\into} & X_T  & \stackrel{f_T}{\to} & T\\
g_U\downarrow\hphantom{g_U}  && g_X\downarrow \hphantom{g_X}&& \hphantom{g}\downarrow g\\
U & \stackrel{j}{\into} & X  & \stackrel{f}{\to} & S\\
\end{array}
\eqno{(\ref{comm.base.ch.defn}.1)}
$$
and a natural base-change  map
$$
g_X^*\bigl(j_*F\bigr)\to  (j_T)_*\bigl(g_U^*F\bigr).
\eqno{(\ref{comm.base.ch.defn}.2)}
$$
We say that $j_*F$ {\it commutes with base change} if
(\ref{comm.base.ch.defn}.2) is an isomorphism for every $g:T\to S$. 

Let  $f:(X,\Delta)\to S$  be a locally stable morphism as in 
Theorem~\ref{half.pgsh.inv.thm}. By  \cite[4.26]{k-modbook},
 there is   a closed subset $Z\subset X$ such that $Z\cap X_s$ has codimension $\geq 2$ in $X_s$ for every $s\in S$
 and both $K_{X/S}$ and 
$\rdown{m\Delta}-D$ are Cartier on $X\setminus Z$. 
In particular, 
$$
\omega_{X/S}^{[m]}\bigl(\rdown{m\Delta}-D\bigr)=
j_* \bigl(\omega_{U/S}^{[m]}(\rdown{m\Delta|_U}-D|_U)\bigr).
$$
Thus we say that
 $\omega_{X/S}^{[m]}\bigl(\rdown{m\Delta}-D\bigr)$  {\it commutes with base change} if $j_* \bigl(\omega_{U/S}^{[m]}(\rdown{m\Delta|_U}-D|_U)\bigr)$
commutes with base change.
\end{defn}

\begin{lem}  \label{comm.base.ch.lem}  Let $(t, T)$ is the spectrum of a DVR,
$f:X\to T$ a finite type morphism of pure relative dimension $n$ and    $j:U\into X$ as in
(\ref{comm.base.ch.defn}). Let $F$ be a coherent sheaf on $U$ that is flat over $T$.  
The following are equivalent.
\begin{enumerate}
\item $j_*F$ commutes with arbitrary base change.
\item $j_*F$ commutes with base change  to $i:\{t\}\into T$.
\item $\depth_{Z_t}(j_*F)\geq 3$.
\end{enumerate}
\end{lem}

Proof.  It is clear that (1) $\Rightarrow$ (2).  The base-change  map
$$
r_t: (j_*F)|_{X_t}\to  (j_t)_*\bigl(F|_{U_t}\bigr).
$$
is an injection and an isomorphism over $U_t$. Furthermore $\depth _{Z_t}(j_t)_*\bigl(F|_{U_t}\bigr)\geq 2$.
Thus $\depth_{Z_t}(j_*F)\geq 3$ iff $r_t$ is surjective. 

Finally  (3) $\Rightarrow$ (1) is straightforward, see   \cite[9.26]{k-modbook}. \qed

\begin{prop} \label{comm.base.ch.prop} Let $S$ be a reduced scheme and  $f:X\to S$ a finite type morphism of pure relative dimension $n$ and    $j:U\into X$ as in
(\ref{comm.base.ch.defn}). 
Let $F$ be a coherent sheaf on $U$ that is flat over $S$. The following are equivalent.
\begin{enumerate}
\item $j_*F$ is flat over $S$ and  commutes with base change.
\item For every morphism  $g:(t,T)\to S$, where $(t, T)$ is the spectrum of a DVR,
$(j_T)_*(g_U^*F\bigr)$  commutes with  base change to $\{t\}\into T$.
\end{enumerate}
\end{prop}
 
Proof. It is clear that (1) $\Rightarrow$ (2). To see the converse,
we use the theory of {\it hulls and husks;} see
\cite{k-hh} or \cite[Chap.9]{k-modbook}.
In this language,   our claim is equivalent to saying that  $j_*F$ is its own hull; cf.\ \cite[9.17]{k-modbook}. 

If $X$ is projective over $S$, then \cite{k-hh} shows that
the hull of $j_*F$ is represented by a monomorphism  $j:S^H\to S$, see also  \cite[9.59]{k-modbook}. Assumption (2)  and Lemma~\ref{comm.base.ch.lem} imply that
whenever $T$ is the spectrum of a DVR and $T\to S$ a morphism whose generic point maps to a generic point of $S$ then $T\to S$ factors through $S^H$. Thus $S^H\to S$ is an isomorphism and so  $j_*F$ is its own hull; cf.\ \cite[3.49]{k-modbook}. 

A similar argument works in general. Assume to the contrary that
$j_*F$ is not its own hull. 
Then there is a point $x\in X$ such that $j_*F$ is not its own hull at $x$.
After completing $X$ at $x$ and $S$ at $f(x)$, we get
$\hat f:\hat X\to \hat S$ such that  
$\hat{j}_*\hat F$ is not its own hull.
Over a complete, local  base scheme the hull functor  is represented by a monomorphism  $j:\hat S^H\to \hat S$ for local maps by \cite[9.61]{k-modbook}. We can now argue as before to get a contradiction. \qed

\section{Small modifications of lc pairs}\label{sec.1}

The local class group $\cl(x\in X)$ of a klt pair $(x\in X, \Delta)$ is finitely generated, and the vector space  $\q\otimes_{\z} \cl(x\in X)$
has a finite  chamber  decomposition into cones such that the small modifications of
$(x\in X)$ are in one-to-one correspondence with the chambers.
The maximal dimensional chambers correspond to the $\q$-factorial small modifications. These were first observed in \cite{r-mmc3} for canonical 3--folds.

For an lc pair $(x\in X)$, the local class group $\cl(x\in X)$ need not be finitely generated. Even if it is and the analogous chamber  decomposition into cones seems to exist, not every chamber  corresponds to a small modification. Already cones over logCY surfaces exhibit many different patterns, as shown by Example~\ref{sm.mod.CY.cnes.exmp}.

Here we study the existence of small modifications of lc pairs. For our applications we need to consider potentially lc pairs as well.

\begin{defn} \label{pot.lc.defn}
Let  $X$ be a normal variety and $\Delta$ an effective $\r$-divisor on $X$. Following \cite{k-db} we say that
$(X, \Delta)$ is {\it potentially lc} if there is an open cover
$X=\cup_i X_i$ and   effective $\r$-divisors  $\Theta_i$  on $X_i$
such that  $(X_i, \Delta|_{X_i}+\Theta_i)$ is lc for every $i$.
Using that most lc surface singularities are rational, hence $\q$-factorial  (cf.\ \cite[10.4 and 10.9]{kk-singbook}) we get the following.
\medskip

{\it Claim \ref{pot.lc.defn}.1.} Let $(X, \Delta)$ be   potentially lc. Then there is a closed subset  $Z\subset X$ of codimension $\geq 3$ such that  $(X\setminus Z, \Delta_{X\setminus Z})$ is lc. \qed
\medskip 

Potentially lc pairs appear quite frequently as auxiliary objects.
For example, if $(Y, \Delta_Y)$ is lc and $\pi:Y\to X$ is a
$(K_Y+ \Delta_Y)$-negative birational contraction then
$(X, \pi_* \Delta_Y)$ is potentially lc.

Canonical models of   non-general type lc pairs have a natural  potentially lc structure, cf.\   \cite{kaw-adj, k-adj}.

(The definition can be naturally generalized to potentially slc pairs, but then (\ref{pot.lc.defn}.1) is not automatic. In order to get a good notion, 
 one probaly should impose
(\ref{pot.lc.defn}.1) as an extra condition.)
\end{defn}

\begin{defn} \label{pot.lc.mod.defn} (cf.\ \cite[1.32]{kk-singbook})
Let  $X$ be a normal variety and $\Delta$ an effective $\r$-divisor on $X$. 
An {\it lc modification} of $(X, \Delta)$ is a proper, birational morphism
$\pi: (X^c, \Delta^c+E^c)\to (X, \Delta)$ where $\Delta^c:=\pi^{-1}_*\Delta$,
$E^c$ is the reduced $\pi$-exceptional divisor, $(X^c, \Delta^c+E^c) $ is lc and
$K_{X^c}+ \Delta^c+E^c$ is $\pi$-ample. 

An  lc modification is unique. As for its existence, we clearly need to assume that $\coeff \Delta\subset[0,1]$. Conjecturally, this is the only necessary condition. 
\cite{MR2955764} shows that  lc modifications exist if $K_X+\Delta$ is $\q$-Cartier. Next we consider some cases when $K_X+\Delta$ is not $\q$-Cartier.
\end{defn}

\begin{prop}\label{pot.lc.can.mod}
 Let $(X, \Delta)$ be a  potentially lc pair. 
Then 
\begin{enumerate}
\item  it has a projective, small, lc modification
$\pi: (X^c, \Delta^c)\to (X, \Delta)$,
\item  $\pi$ is a local isomorphism at every lc center of $(X^c, \Delta^c) $ and
\item $\pi$ is a local isomorphism over $x\in X$ iff $K_X+\Delta$ is $\r$-Cartier at $x$. 
\end{enumerate}
\end{prop}

Proof. Since an lc  modification is unique, it is enough to construct it locally on $X$. We may thus assume that there is an  effective $\r$-divisor $\Theta$  on $X$
such that  $(X, \Delta+\Theta)$ is lc. 

By  \cite{k-db, fujino-ssmmp},  $(X, \Delta+\Theta)$ has a $\q$-factorial dlt modification
$(X', \Delta'+\Theta'+E')\to (X, \Delta+\Theta)$.

Next we have a canonical model for 
$(X', \Delta'+E')\to (X, \Delta)$ by \cite[1.6]{MR3032329}, call it
$\pi: (X^c, \Delta^c+E^c)\to (X, \Delta)$.
Thus $-\Theta^c\sim_{\r,\pi}K_{X^c}+\Delta^c+E^c$ is $\pi$-ample. 
In particular, $\ex(\pi)\subset \supp \Theta^c$  by Lemma~\ref{neg.exc.lem}. Thus $\pi$ is small,
$E^c=0$ and $\pi: (X^c, \Delta^c)\to (X, \Delta)$ is the required  lc modification. 

Let $W\subset X^c$ be an lc center of $(X^c, \Delta^c)$.  Increasing
$ \Delta^c$ to $\Delta^c+ \Theta^c$ decreases discrepancies, and the decrease is strict for divisors whose center is contained in $\supp \Theta^c$.
Since $(X^c, \Delta^c+ \Theta^c)$ is lc, this implies that
$W\not\subset \supp \Theta^c$. Since $\ex(\pi)\subset \supp \Theta^c$, this shows that $W\not\subset\ex(\pi)$. 

If $K_X+\Delta$ is $\r$-Cartier at $x$ then $K_{X^c}+ \Delta^c\simr \pi^*(K_X+\Delta)$ near $x$. Since $K_{X^c}+ \Delta^c$ is  $\pi$-ample,
this implies that $\pi$ is a local isomorphism over $x$. 
\qed

\begin{lem}\label{neg.exc.lem} Let $\pi:Y\to X$ be a proper birational morphism, $X$ normal.
Let $D$ be an effective divisor on $X$ such that $-\pi^{-1}_*D$ is $\pi$-nef. Then $\supp\pi^{-1}D= \supp\pi^{-1}_*D$.
If $-\pi^{-1}_*D$ is $\pi$-ample then $\ex(\pi)\subset \supp \pi^{-1}_*D$ and $\pi$ is small. \qed
\end{lem}
\medskip

As we see in (\ref{sm.mod.CY.cnes.exmp}.1), 
 an lc pair usually does not have a small, $\q$-factorial modification. 
However, we can at least achieve that all irreducible components of $\Delta$ become $\q$-Cartier.

\begin{cor}  \label{lc.qf.cor.1}
Let $(X, \Delta)$ be an lc pair. Then  $(X, \Delta)$ has a projective, small modification $\pi: \bigl(X^w, \Delta^w:=\pi^{-1}_*(\Delta)\bigr)\to (X, \Delta)$ such that
\begin{enumerate}
\item $(X^w, \Delta^w)$ is lc,
\item $K_{X^w}$ is $\pi$-nef, 
\item every irreducible component of $\Delta^w$ is $\q$-Cartier,
\item $\ex(\pi)\subset \supp \Delta^w$ and 
\item  $\pi$ is an isomorphism over $x\in X$ iff every irreducible component of $\Delta$ is $\q$-Cartier at $x$.
\end{enumerate}
\end{cor}

Proof.  If $\pi$ is small then $K_{X^w}+\Delta^w\simr \pi^*(K_X+\Delta)$ hence 
$(X^w, \Delta^w)$ is lc.

First  apply
Proposition~\ref{pot.lc.can.mod} to  $(X, 0)$ to get  $\tau_1: X_1\to X$ 
such that $K_{X_1}$ is $\tau$-ample and $\ex(\tau)$ is contained in the support of $ \Delta_1:=(\tau_1)^{-1}_*\Delta$. 
If $\tau_2:X_2\to X_1$ is any small modification then
$K_{X_2}\simr \tau_2^*K_{X_1}$, thus $K_{X_2}$ is $\tau_1\circ \tau_2$-nef.
Similarly, since $\Delta_1$ is $\r$-Cartier, $\Delta_2:=(\tau_2)^{-1}_*\Delta_1\simr \tau_2^*\Delta_1$, hence (4) also  holds for us. 

In order to avoid multi-indices, we can thus assume that $K_X$ is $\q$-Cartier.
 
We need  to achieve (3).
Let $a_1D_1$ be an  irreducible component of $\Delta$ and  apply
Proposition~\ref{pot.lc.can.mod} to  $(X, \Delta-a_1D_1)$ to get  $\sigma: (X', \Delta'-a_1D'_1)\to (X, \Delta-a_1D_1)$. Since  $K_{X'}+\Delta'-a_1D'_1$
and $K_{X'}+\Delta'\simr \sigma^*(K_{X}+\Delta)$ are $\r$-Cartier, so is
$D'_1$ and $\ex(\sigma)\subset \supp D'_1$. 

By induction on the
number of irreducible components of $\Delta$ we may assume that
the claims (1--5) hold for $(X', \Delta'-a_1D'_1)$ and we have
$$
\pi': \bigl(X^w, \Delta^w-a_1D^w_1\bigr)\to (X', \Delta'-a_1D'_1).
$$
 The composite
$$
\pi:=\sigma\circ \pi': \bigl(X^w, \Delta^w\bigr)\to (X, \Delta)
$$
is the required modification. \qed


\medskip

Note that $ \bigl(X^w, \Delta^w\bigr)$  is not unique and the different
ones are related to each other by flops. However, as shown by \cite[96]{k-exrc} and (\ref{sm.mod.CY.cnes.exmp}.3),  flops do not always exist in the lc case, so we have less freedom than in the klt case.

Proposition~\ref{pot.lc.can.mod} does not extend to
potentially slc pairs. In codimension 2 this is related to the bookkeeping problem we already encountered in  Paragraph~\ref{half.rest.well.lem}. For example, the pair 
$(S,D):=\bigl((xy=0), (x=z=0)\bigr)$ is potentially slc since
$\bigl((xy=0), (z=0)\bigr)$ is slc, but $K_S+D$ is not $\q$-Cartier.
Since a demi-normal surface does not have small modifications, there is nothing that can be done.
However, even if our divisor is $\q$-Cartier in  codimension 2, there are higher codimension obstructions. 
The following  is a slight modification of the second  example in \cite[1.40]{kk-singbook}. 

\begin{exmp}\label{1.40.kk-singbook.exmp} 
  Let $Z\subset \p^{n-1}$ be a smooth hypersurface of degree $n$ and $D_1\subset
  X_1:=\a^n$ the cone over $Z$.  
Next fix two points $p,q\in \p^1$, embed $Z\times \p^1$ into $\p^{2n-1}$ by the global
  sections of $\o_{Z\times \p^1}(1,1)$ and let $\bigl(X_2, D_2+D'_2\bigr)$ be the cone
  over $\bigl((Z\times \p^1), (Z\times \{p\})+(Z\times \{q\}) \bigr)$.

 $\bigl(X_1, D_1\bigr)$ and $\bigl(X_2, D_2+D'_2\bigr)$ are both lc by
\cite[3.1]{kk-singbook}, thus 
  one can glue them using the natural
  isomorphism $\sigma: D_1\cong D_2$ to obtain an slc pair 
$\bigl(X_1\amalg_{\sigma}X_2, D'_2\bigr)$. 
Thus $(X, 0):=\bigl(X_1\amalg_{\sigma}X_2, 0\bigr)$ is potentially slc. 

We check in \cite[1.40]{kk-singbook} that $(X, 0)$
  has no slc modification. The reason is that 
$\bigl(X_1, D_1\bigr)$ is lc, hence its lc modification is itself, but
the lc modification of $\bigl(X_2, D_2\bigr)$ is obtained by first taking the conical resolution and then  contracting the
  $\p^1$-factor of the exceptional divisor  $E\cong Z\times p^1$.
Thus the birational transform of $D_1$ in the lc modification is
isomorphic to $D_1$ but the birational transform of $D_2$ in the lc modification is
isomorphic to its blow-up. So the two lc modifications can not be glued together.

A similar example, where the canonical class is Cartier is the following.

Let $v\in Q\subset \a^4$ be a quadric cone with vertex  $v$ and 
$v\in H\subset Q$ a hyperplane section. Assume that both  $Q$ and $H$ have an isolated singularity at $v$. We can glue  
$(Q, H)$ and  $(H\times \a^1, H\times \{0\})$ to get an semi-dlt 3-fold $X$
(cf.\ \cite[5.19]{kk-singbook})
and $K_X\sim 0$. 

Next we add a boundary. Let
 $A\subset Q$ be a plane. 
Then $L:=A\cap H$ is line through $v$. The boundary $\Delta$ is $A$ on $Q$ and
$L\times \a^1$ on $H\times \a^1$. 
Note that $L\times \a^1$ is $\q$-Cartier on $H\times \a^1$ but $A$ is not
$\q$-Cartier on $Q$. The lc modification of $(Q, H+A)$ is one of the
small modifications of $Q$ and the birational transform of $H$ is its minimal resolution. Thus again the 2 small modifications can not be glued together.

\end{exmp}

The next result shows that for locally stable morphisms, the existence of an 
slc modification depends only on the generic fiber.

\begin{prop} Let $f:(X, \Delta+\Delta')\to B$ be a locally stable morphism to a smooth, irreducible curve $B$.
Assume that the generic fiber   has an slc modification $\pi_g: (X_g^c, \Delta_g^c)\to (X_g, \Delta_g)$.
Then $(X, \Delta)$ also has an slc modification.
\end{prop}

Proof. We follow the method of \cite[Sec.2.4]{k-modbook}. That is, first we normalize $(X, \Delta)$, then use Proposition~\ref{pot.lc.can.mod} to obtain the
lc modification of the normalization and finally use the gluing theory of
\cite{k-source} and \cite[Chap.9]{kk-singbook} to get  $\pi: (X^c, \Delta^c)\to (X, \Delta)$.

Thus we start with  the normalization
$$
\rho: (\bar X, \bar\Delta+\bar D+\bar\Delta')\to (X, \Delta+\Delta').
$$
By Proposition~\ref{pot.lc.can.mod}, 
there is a  projective, small, lc modification
$$
\sigma: (\bar X^c, \bar\Delta^c+\bar D^c)\to (\bar X, \bar\Delta+\bar D),
$$
whose generic fiber over $B$ is the normalization of
$(X_g^c, \Delta_g^c)$. Let  $\bar D^{cn}\to \bar D^c$ denote the normalization.
On the generic fiber we have an involution
$$
\tau_g:  \bigl(\bar D^{cn}_g, \diff \bar\Delta^c_g\bigr)\leftrightarrow
\bigl(\bar  D^{cn}_g, \diff \bar\Delta^c_g\bigr),
$$
where we take the different on $\bar D^{cn}_g $.
By \cite[2.12]{k-modbook} $\bigl(\bar D^{cn}_g, \diff \bar\Delta^c\bigr) $ has no log centers supported over a closed point of $B$, thus,  
by \cite[2.45]{k-modbook}, $\tau_g$   extends to an involution
$$
\tau:  \bigl(\bar D^{cn}, \diff \bar\Delta^c\bigr)\leftrightarrow
\bigl(\bar D^{cn}, \diff \bar\Delta^c\bigr),
$$
where we take the different on $\bar D^{cn} $.

By \cite[2.54]{k-modbook} 
we thus obtain  $(X^c, \Delta^c)$ as the geometric quotient of
$(\bar X^c, \bar\Delta^c+\bar D^c) $ by the equivalence relation generated by $\tau$.  \qed

\section{Proof of Propositions  \ref{tech.prop.0} and \ref{tech.prop.B}}\label{sec.2}

We prove Proposition~\ref{tech.prop.2}, which is a generalization of   Proposition~\ref{tech.prop.0}. Its assumptions are somewhat convoluted, but probably sharp; see Examples \ref{half.main.exmp}.6--7. They were arrived at by looking at a proof of  Proposition~\ref{tech.prop.0} and then trying to write down a minimal set of assumptions that make the arguments work. For now the only applications I know of are Propositions~\ref{tech.prop.0}--\ref{tech.prop.B} and  Theorem~\ref{half.pgsh.inv.thm}. 

We start with the following application of  the method of Section \ref{sec.1}. 

\begin{prop} \label{tech.prop.1}
Let $(X,\Delta)$ be an lc pair and $\Delta_1, \Delta_2$ effective divisors such that  $\Delta_1+ \Delta_2\leq \Delta$. Let $B$ be a Weil $\z$-divisor such that  $B\simr K_X+L+\Delta_1-c\Delta_2$ where $L$ is $\r$-Cartier and $c\geq 0$. Then there is a small, lc modification
$\pi:(X', \Delta')\to (X,\Delta)$  such that
\begin{enumerate}
\item  $B'$ is  $\q$-Cartier,
\item $K_{X'}+\Delta'_1$ is $\r$-Cartier,
\item  $\ex(\pi)\subset \supp (\Delta'-\Delta'_1)$,
\item none of the lc centers of $(X',\Delta'_1)$ are contained in $\ex(\pi)$, 
\item $-\Delta'_2$ is $\r$-Cartier and $\pi$-nef,
\item $R^i\pi_*\o_{X'}(B')=0$ for $i>0$ and
\item  $H^i\bigl(X, \o_X(B)\bigr)=H^i\bigl(X,\o_{X'}(B')\bigr)$.
\end{enumerate}
\end{prop}

Proof.  We construct $\pi:(X', \Delta')\to (X,\Delta)$ in 2 steps.
First we apply Proposition~\ref{pot.lc.can.mod} to   $(X, \Delta-\Delta_2)$.
We get  $\tau_1:(X^*, \Delta^*)\to (X,\Delta)$ such that
$$
-\Delta_2^*\simr K_{X^*}+\Delta^*-\Delta^*_2- \tau_1^*(K_{X}+\Delta)
$$ 
is $\tau_1$-ample and its support contains the exceptional set of $\tau_1$. 
Since $\Delta_1^*\leq \Delta^*-\Delta_2^*$, we can 
next  apply Proposition~\ref{pot.lc.can.mod} to   $(X^*, \Delta_1^*)$
to get $\tau_2: (X', \Delta')\to (X^*, \Delta^*)$. Set $\pi:=\tau_1\circ \tau_2$.
By construction  $\Delta'_2=\tau_2^*\Delta_2^*$ and $K_{X'}+\Delta'_1$ are $\r$-Cartier, and so is 
$B'\simr (K_{X'}+\Delta'_1)+\pi^*L-c\Delta'_2$.

Furthermore,  $\ex(\pi)$ is contained in
$\supp (\Delta'_2)\cup \supp (\Delta'-\Delta'_1)$.
Since $\Delta'_2\leq \Delta'-\Delta'_1$, this implies that
 $\ex(\pi)\subset \supp (\Delta'-\Delta'_1)$. Thus
none of the lc centers of $(X',\Delta'_1)$ are contained in $\ex(\pi)$.
Since $-\Delta'_2=-\tau_2^*\Delta_2^*$ is $\pi$-nef, so is
$\pi^*L-\Delta'_2$ 
and these in turn imply that  $R^i\pi_*\o_{X'}(B')=0$ for $i>0$ by Theorem~\ref{fujino.1.10.thm}. 
Finally the Leray spectral sequence shows (7). \qed
\medskip

It is convenient to state the following results using the notion of 
Mumford divisors; see \cite{k-mumf} for a detailed treatment. 

\begin{defn} \label{mumf.div.defn}
Let $X$ be a scheme.  I call a  Weil divisor $B$ on $X$ a
{\it Mumford divisor} if $X$ is regular at all generic points of $\supp B$.
Thus on a normal variety, every  Weil divisor is a Mumford divisor.



Let $S\subset X$  be a closed subscheme.
 $B$ is called a {\it Mumford divisor  along $S$} if 
 $\supp B$ does not contain any 
irreducible component of $S$,
  $B$ is Cartier  at all generic points of $S\cap \supp B$ and
 $S$ is regular at all generic points of $S\cap \supp B$.

These imply that $B|_S$ is a well-defined Mumford divisor on $S$ and 
there is a subset $Z\subset S$ of codimension $\geq 2$ such that 
$B$ is Cartier at all points of  $S\setminus Z$. The restriction
 sequence
$$
0\to \o_X(B)(-S)\to \o_X(B)\to \o_S(B|_S)\to 0,
\eqno{(\ref{mumf.div.defn}.1)}
$$
 is left exact everywhere and right exact on $X\setminus Z$. 
(We use this sequence only when $X$ and $S$ are both $S_2$; then there is no ambiguity in the definition of $\o_X(B)$ and $\o_S(B|_S)$.)
Assume next that $S\subset X$ is a Cartier divisor. Then
(\ref{mumf.div.defn}.1) gives a natural injection
$$ 
r: \o_X(B)|_S\into \o_S(B|_S),
\eqno{(\ref{mumf.div.defn}.2)}
$$
which is an isomorphism on $S\setminus Z$. Since $\o_S(B|_S) $ is $S_2$ by definition, $r$ is  an isomorphism everywhere
 iff $\o_X(B)|_S$ is $S_2$. As we noted in Definition~\ref{S3.on.S.defn}, this is equivalent to
$\o_X(B)$ being $S_3$ along $S$. We thus obtain the following observation.
\medskip

{\it Claim \ref{mumf.div.defn}.3.} If $S\subset X$ is a Cartier divisor then the sequence (\ref{mumf.div.defn}.1) is exact iff $\o_X(B)$ is $S_3$ along $S$.\qed
\end{defn}

We need to understand the $S_3$ condition for divisorial sheaves on slc pairs. 
The  first part  of the 
following is proved  in \cite{ale-lim} and  the second part in \cite{k-dep}; see also
\cite[7.20]{kk-singbook} and \cite[Sec.7.1]{fuj-book}. 

\begin{thm}\label{alex.kol.depth.3} Let $(X,\Delta)$ be an slc pair and $x\in X$ a point that is not an lc center. Let $B$ be a Mumford $\z$-divisor on $X$. Assume that 
\begin{enumerate}
\item either $B$ is $\q$-Cartier,
\item or  $B\simr -\Delta'$ for some $0\leq \Delta'\leq \Delta$.
\end{enumerate}
Then
 $\depth_x \o_X(B)\geq \min\{3, \codim_Xx\}$. \qed
\end{thm}

\begin{cor}\label{alex.kol.depth.3.cor} Let $(X,\Delta)$ be an slc pair,
$S$ a reduced $\q$-Cartier, Mumford  divisor and $B$ a $\q$-Cartier $\z$-divisor that is Mumford along $S$. Assume that 
$B$ satisfies (\ref{alex.kol.depth.3}.1) or (\ref{alex.kol.depth.3}.2) and 
$\supp S\subset \supp \Delta$. 
Then the sequence
$$
0\to \o_X(B-S)\to \o_X(B)\to \o_S(B|_S)\to 0 \qtq{is exact.}
$$
\end{cor}

Proof. Assume first that $S$ is Cartier.  $(X,\Delta-\epsilon S)$ is also an slc pair for $0<\epsilon\ll 1$ and
none of its lc centers are contained in $S$ by \cite[2.27]{km-book}. Thus
$\o_X(B)$ is $S_3$ along $S$ by (\ref{alex.kol.depth.3}.1), hence
the  sequence is exact by (\ref{mumf.div.defn}.3). 

If $S$ is not Cartier, let $m>0$ be the smallest integer such that $mS\sim 0$. This linear equivalence defines a degree $m$ cyclic cover
$\tau:\bar X\to X$ such that $\bar S:=\tau^*(S)$ is Cartier; see \cite[2.49--53]{km-book}. Set $\bar B:=\tau^*(B)$. We have already established that  $\depth_{\bar x} \o_{\bar X}(\bar B)\geq \min\{3, \codim_{\bar X}\bar x\}$ for every $\bar x\in \bar S$.
Since  $\o_X(B)$ is a direct summand of 
$  \tau_*\o_{\bar X}(\bar B)$, this implies that $\o_X(B)$ is $S_3$ along $S$.  \qed

\medskip

Now we come to the main technical result, which is a strengthening of
Proposition~\ref{tech.prop.B}. Examples \ref{half.main.exmp}.6--7
shows that the assumptions are likely optimal.

\begin{prop} \label{tech.prop.2}
Let $(X,S+\Delta)$ be an lc pair where $S$ is $\q$-Cartier.  Let $B$ be a Weil $\z$-divisor that is Mumford along $S$ and $\Delta_3$ an effective $\r$-divisor such that 
\begin{enumerate}
\item $B\simr -\Delta_3$,
\item $\Delta_3\leq \Delta^{(=1/2)}+ \rup{\Delta^{(> 1/2)}}$ and 
\item $\rdown{\Delta_3}\leq \rdown{\Delta}$.
\end{enumerate}
Then  $\o_X(B)$ is $S_3$ along $S$.
\end{prop}

Proof. Assume first that  $S$ is Cartier with equation $s=0$.
By (\ref{mumf.div.defn}.3) we need to show that the sequence (\ref{mumf.div.defn}.1) is exact.

By Theorem~\ref{alex.kol.depth.3} 
  we need to focus on the points where $B$ is not $\q$-Cartier. This suggests that we should use Proposition~\ref{tech.prop.1}. The question is local on $X$, so we may  assume that $K_X+\Delta\simr 0$. 
By assumption (1)
$$
B\simr -\Delta_3\simr K_X+\Delta -\Delta_3.
\eqno{(\ref{tech.prop.2}.4)}
$$ 
Set $\Delta_1:=\bigl(\Delta -\Delta_3\bigr)^{\geq 0}$ and
$\Delta_2:=\epsilon\bigl(\Delta -\Delta_3\bigr)^{< 0}$
for some $0<\epsilon\ll 1$. 
 It is clear that 
 $\supp \Delta_1$, $\supp \Delta_2$ have no common irreducible components
and $\Delta_1+ \Delta_2\leq \Delta$ for $0<\epsilon\ll 1$.
Furthermore, $B\simr K_X+\Delta_1 -c\Delta_2$ with $c:=\epsilon^{-1}$.
Using these $\Delta_1, \Delta_2$ in 
 Proposition~\ref{tech.prop.1}
 we  obtain $\pi:(X', \Delta')\to (X,\Delta)$.

Note that  $B'$ is $\q$-Cartier by  (\ref{tech.prop.1}.1),
  $(X', S'+\Delta')$ is lc and
$\ex(\pi)\subset \supp (\Delta'-\Delta'_1)$ by (\ref{tech.prop.1}.3).
Thus none of the lc centers of $(X', S'+\Delta'_1)$ are contained in $\ex(\pi)$, in particular, $S'$ is smooth at the generic points of all exceptional divisors of $\pi_S:=\pi|_{S'}:S'\to S$. Thus  $B'$ is also Mumford along $S'$, hence  the sequence
$$
0\to \o_{X'}(B')\stackrel{s}{\to} \o_{X'}(B')\to \o_{S'}\bigl(B'|_{S'}\bigr)\to 0
\eqno{(\ref{tech.prop.2}.5)}
$$
is exact by Theorem~\ref{alex.kol.depth.3.cor}. Since   $R^1\pi_*\o_{X'}(B')=0$  by (\ref{tech.prop.1}.6), 
pushing (\ref{tech.prop.2}.5) forward gives an exact  sequence
$$
0\to \pi_*\o_{X'}(B')\stackrel{s}{\to} \pi_*\o_{X'}(B')\to (\pi_{S})_*\o_{S'}\bigl(B'|_{S'}\bigr)\to 0.
\eqno{(\ref{tech.prop.2}.6)}
$$
It remains  to show that  
$(\pi_S)_*\o_{S'}\bigl(B'|_{S'}\bigr)=\o_{S}\bigl(B|_{S}\bigr)$. 
This holds if there is a $B''\leq B'|_{S'}$ such that 
$(\pi_S)_*\o_{S'}\bigl(B''\bigr)=\o_{S}\bigl(B|_{S}\bigr)$.

Assume first that $\rdown{\Delta_3}=0$ and
$\Delta_3\leq \rup{\Delta^{(> 1/2)}}$.
By assumption  $B+\Delta_3\simr 0$, hence $B'+\Delta'_3\simr 0$.
Thus $(B'+\Delta'_3)|_{S'}\simr 0$ and using Lemma~\ref{push.wfd.lem}
with $N:=-(B'+\Delta'_3)|_{S'}$ and $H:=0$ gives that 
$$
(\pi_S)_*\o_{S'}\bigl(\rdown{B'|_{S'}+\Delta'_3|_{S'}}\bigr)=\o_S\bigl(\rdown{B|_{S}+\Delta_3|_{S}}\bigr).
\eqno{(\ref{tech.prop.2}.7)}
$$
Since $S'$ is Cartier, $B'|_{S'}$ is a $\z$-divisor, so
$$
\rdown{B'|_{S'}+\Delta'_3|_{S'}}=B'|_{S'}+\rdown{\Delta'_3|_{S'}}
\qtq{and}
\rdown{B|_{S}+\Delta_3|_{S}}=B|_{S}+\rdown{\Delta_3|_{S}}.
$$
Since we assume that $\rdown{\Delta_3}=0$,
   (\ref{half.rest.well.lem}.3) and our assumption 
$\Delta_3\leq \rup{\Delta^{(> 1/2)}}$  imply that
$\rdown{\Delta_3|_{S}}=\rdown{\Delta_3}|_{S}=0$ and 
$\rdown{\Delta'_3|_{S'}}=\rdown{\Delta'_3}|_{S'}=0$.
Thus
$$
(\pi_S)_*\o_{S'}\bigl(B'|_{S'}\bigr)=\o_S\bigl(B|_{S}\bigr)
\eqno{(\ref{tech.prop.2}.8)}
$$
and we are done with this case.

In general, let $E\subset S'$ be the largest, reduced, $\pi_s$-exceptional divisor
such that $\pi_S(\supp E)\subset \supp \Delta_3$.  
Set $H:=\rdown{\Delta_3|_S}$ and let $H'$ denote the birational transform of 
$H$ on $X'$. 

We apply 
Lemma~\ref{push.wfd.lem}
with $N:=-(B'+\Delta'_3)|_{S'}$ and $H'$ to get that 
$$
(\pi_S)_*\o_{S'}\bigl(\rdown{B'|_{S'}+\Delta'_3|_{S'}-H'-\epsilon E}\bigr)=\o_S\bigl(\rdown{B|_{S}+\Delta_3|_{S}-H}\bigr).
\eqno{(\ref{tech.prop.2}.9)}
$$
This in turn gives (\ref{tech.prop.2}.8) provided that we prove that
$$
\rdown{B'|_{S'}+\Delta'_3|_{S'}-H'-\epsilon E}\leq B'|_{S'}\qtq{and}
\rdown{B|_{S}+\Delta_3|_{S}-H}=B|_{S}.
\eqno{(\ref{tech.prop.2}.10)}
$$
The second of these holds by our choice of $H$  since
$$
B|_{S}+\Delta_3|_{S}-H=B|_{S}+\Delta_3|_{S}-\rdown{\Delta_3|_S}=
B|_{S}+\{\Delta_3|_{S}\}.
\eqno{(\ref{tech.prop.2}.11)}
$$
The same argument shows that $\rdown{\Delta'_3|_{S'}-H'}$ is $\pi_S$-exceptional. It remains to understand what happens along any $\pi_S$-exceptional prime divisor $F\subset S'$. 
If $F\subset \supp \Delta'_3$ then  also 
$F\subset \supp \Delta'$ and,  
as we argue in the proof of (\ref{half.rest.well.lem}.4),
 we are in one of the following cases.
\begin{enumerate}
\item[(a)]  $\coeff_F\Delta'_3|_{S'}<1$,
\item[(b)]  $F$ is contained in a unique
irreducible component  $D'_i$ of $\Delta'$,   $S'\cap D'_i$ has multiplicity 1 along $F$ and $\coeff_{D_i}\Delta_3=\coeff_{D_i}\Delta=1$,
\item[(c)]     $F$ is contained in a unique
irreducible component  $ D'_i$ of $\Delta'$,  $S'\cap D'_i$ has multiplicity 2 along $F$  and $\coeff_{D_i}\Delta_3=\coeff_{D_i}\Delta=\frac12$, 
\item[(d)]     $F$ is contained in  two
irreducible components  $D'_i,D'_j $ of $\Delta'$, $S'\cap D'_i, S'\cap D'_j$ both have  multiplicity 1 along $F$   and $\coeff_{D_i}\Delta_3=\coeff_{D_i}\Delta=\coeff_{D_j}\Delta_3=\coeff_{D_j}\Delta=\frac12$. 
\end{enumerate}

In the first  case 
$F\not \subset \rdown{\Delta'_3|_{S'}}$. 

If (b) (resp.\ (c) or (d)) holds then $\coeff_F(\Delta'_3|_{S'})=1$ by assumption (3) (resp.\ (2)). 
Moreover,   $\pi_S(\supp F)\subset \supp(D_i|_S)$  and
$D_i|_S$ appears in $\Delta_3|_S$ with coefficient $1$ or $\frac12$. 
Thus  $\pi_S(\supp F)\subset \supp\bigl(\pi_*\{-(B'+\Delta'_3)|_{S'}\}+\pi_*H'\bigr)$, hence 
$F\subset E$. Thus 
$$
\coeff_F\bigl(\Delta'_3|_{S'}-H'-\epsilon E\bigr)=1-\epsilon.
$$
Using this for every $F$ we get that 
$$
\rdown{B'|_{S'}+\Delta'_3|_{S'}-H'-\epsilon E}=B'|_{S'},
$$
hence  (\ref{tech.prop.2}.8) holds. 
This completes the proof if $S$ is Cartier.




If $S$ is only $\q$-Cartier, a suitable cyclic cover reduces everything to the Cartier case, as in the proof of Corollary~\ref{alex.kol.depth.3.cor}. \qed

\medskip
The following is a related to the negativity lemma \cite[3.39]{km-book}.

\begin{lem} \label{push.wfd.lem}
Let $\pi:Y\to X$ be a proper, birational morphism of normal schemes.
Let $N, H$ be $\r$-divisors such that $N$ is $\pi$-nef and
$H$ is effective and horizontal.
Then  $$
\pi_*\o_Y(\rdown{-N-H})=\o_X(\rdown{\pi_*(-N-H)}).\eqno{(\ref{push.wfd.lem}.1)}
$$
 Furthermore, let $E$ be the largest reduced,  effective $\pi$-exceptional divisor such that
$\pi(\supp E)\subset \supp(\pi_*\{N\}+ \pi_* H)$. Then
 $$
\pi_*\o_Y(\rdown{-N-H-\epsilon E}) =\o_X(\rdown{\pi_*(-N-H)})\eqno{(\ref{push.wfd.lem}.2)}
$$
for $0<\epsilon\ll 1$.
\end{lem}

Proof.  The question is local on $X$, so we may assume that $X$ is affine. By the Chow lemma we may assume that $\pi$ is projective. 
Let $F$ denote the divisorial part of $\ex(\pi)$.  

Let  $s$ be a section of $\o_X(\rdown{\pi_*(-N-H)})$.
Then $\pi^*s$ is a section of $\o_Y(\rdown{-N-H}+mF)$ for some $m\geq 0$ and
(\ref{push.wfd.lem}.1) is equivalent to saying that $\pi^*s$ has no poles along $F$.  If this fails then 
we can  cut $Y$ with general hyperplanes, until (after Stein factorization) we get a counter example $\tau:Y'\to X'$ of dimension 2.
In this case  $N':=N|_{Y'}$ and $H':=H|_{Y'}$ are both $\tau$-nef  $\r$-divisors.

Another  advantage of the 2-dimensional normal situation is that there is a well-defined, monotone  pull-back  for all Weil divisors. That is, if $B_1\leq B_2$ then $\tau^*B_1\leq \tau^*B_2$.  Furthermore, there are effective $\tau$-exceptional divisors $E_i$ such that
$$
\begin{array}{l}
\tau^*\tau_*N'=N'+E_1,\quad
\tau^*\tau_*H'=H'+E_2\qtq{and}\\ 
\tau^*\tau_*\bigl(\rup{N'+H'}-N'-H'\bigr)=\bigl(\rup{N'+H'}-N'-H'\bigr)^h+E_3,
\end{array}
\eqno{(\ref{push.wfd.lem}.3)}
$$
where the first claim uses \cite[3.39]{km-book},
 the superscript $h$ denotes the horizontal part and
$\supp E=\supp E_2\cup\supp E_3$. 
Note that 
$$
\rdown{-N'-H'}=-N'-H'-\bigl(\rup{N'+H'}-N'-H'\bigr).
\eqno{(\ref{push.wfd.lem}.4)}
$$
Putting these together gives that 
$$
\begin{array}{rcl}
\tau^*\rdown{\tau_*(-N'-H')}& = & -\tau^*\tau_*N' -\tau^*\tau_*H'-
\tau^*\tau_*\bigl(\rup{N'+H'}-N'-H'\bigr)\\
& = & -N'-H'-\bigl(E_1+E_2+E_3+\bigl(\rup{N'+H'}-N'-H'\bigr)^h\bigr)\\
& \leq & -N'-H' -\epsilon E
\end{array}
$$
for $0<\epsilon \ll 1$. Let now $s'$ be any section of
$\o_{X'}(\rdown{\tau_*(-N'-H')})$. 
 Then we get that 
$$
(\tau^*s')\leq \tau^*\rdown{\tau_*(-N'-H')}\leq -N'-H' -\epsilon E.
\eqno{(\ref{push.wfd.lem}.5)}
$$
Since $(\tau^*s')$ is a $\z$-divisor, we also have the stronger inequality 
$$
(\tau^*s')\leq \rdown{-N'-H' -\epsilon E}.
\eqno{(\ref{push.wfd.lem}.6)}
$$
This implies (\ref{push.wfd.lem}.2). \qed

\medskip

In some situations one may need the following   variants of 
Lemma~\ref{push.wfd.lem}. The first one  can be obtained by the same proof and the second one can be reduced to the normal case  once we  assume that the singularities do not interfere with the divisors much. 
\medskip

{\it Complement \ref{push.wfd.lem}.7.}
Using the above notation and assumptions, let 
 $E^*$ be the divisorial part of $\ex(\pi)$  and  fix a point $x\in X$. Then,  for  $0<\epsilon\ll 1$, we have 
\begin{enumerate}
\item either $\pi_*\o_Y(\rdown{-N-H-\epsilon E^*})=\o_X(\rdown{\pi_*(-N-H)})$
near $x$ 
\item or $N$ is a principal divisor near $\pi^{-1}(x)$  and   $\supp H\cap \pi^{-1}(x)=\emptyset$. \qed
\end{enumerate}

\medskip

{\it Complement \ref{push.wfd.lem}.8.}
Let $\pi:Y\to X$ be a proper, birational morphism of pure dimensional,  reduced schemes. Assume that $Y$ is $S_2$ and   $\pi$ is a local isomorphism at  all  codimension 1 singular points  of $X$ or $Y$.
Let $N$ be a $\pi$-nef Mumford divisor on $Y$ 
and $H$ an effective,  horizontal Mumford $\r$-divisor.
Then  (\ref{push.wfd.lem}.1--2)  hold. \qed
\medskip

\begin{say}[Proof of  Proposition \ref{tech.prop.0}]\label{pf.of.tech.prop.0}
We may assume that $X$ is affine.
Let $x\in H$ be a point of codimension 1. Then either $H$ and $X$ are both smooth at $x$ or $H$ has a node and $X$ a Du~Val singularity at $x$. In the latter case $x\not\in \supp\Delta$. Thus  $mK_X+\rdown{m\Delta}$ is Cartier at $x$, hence a general divisor $B\sim mK_X+\rdown{m\Delta}-D$ is Mumford along $H$.

In order to prove Proposition~\ref{tech.prop.0} we apply  Proposition~\ref{tech.prop.2} to $B$ 
with $\Delta_3:=m\Delta-\rdown{m\Delta}+D$. Thus 
$$
B\sim mK_X+\rdown{m\Delta}-D= m(K_X+\Delta) -\Delta_3\simr  -\Delta_3,
$$
and the proof of Proposition~\ref{tech.prop.2} uses  Proposition~\ref{tech.prop.1} with
$$
\Delta_1:=\bigl(\rdown{m\Delta} -(m-1)\Delta\bigr)^{\geq 0}\qtq{and}
\Delta_2:=\epsilon\bigl(\rdown{m\Delta} -(m-1)\Delta\bigr)^{< 0}.
$$
It is  clear that $\Delta_3\leq \rup{\Delta^{(\geq 1/2)}}$, but 
(\ref{tech.prop.2}.2) needs  a stronger statement for divisors whose
coefficient is $\frac12$. If $A$ is a prime divisor on $X$ such that  
 $\coeff_A\Delta=\frac12$ then $\coeff_A\Delta_3=0$ if  $m$ is even 
and $\coeff_A\Delta_3=\frac12$ if $m$ is odd. Thus, in all cases,
$$
\coeff_A\Delta_3\leq \tfrac12 =\coeff_A\bigl(\Delta^{(=1/2)}+\rup{\Delta^{(>1/2)}}\bigr),
$$
as required for (\ref{tech.prop.2}.2).
So the assumptions of Proposition~\ref{tech.prop.2} are satisfied and
we get Proposition~\ref{tech.prop.0}. \qed
\end{say}




\section{Global applications}\label{sec.3}

The previous constructions  can also  be used to obtain vanishing theorems for
Weil divisors that are not assumed to be $\q$-Cartier.
It turns out that the Kawamata-Viehweg vanishing, even in its strong form given in \cite[1.10]{MR3238112}, holds without the $\q$-Cartier assumptions, see Corollary~\ref{FKV.cor}. The  proof below  uses too much of MMP; it would be desirable to have an argument without such heavy tools.

\begin{defn} Let $(X, \Delta)$ be a  potentially lc pair.
An irreducible subvariety $W\subset X$ is an {\it lc center} of $(X, \Delta)$
if there is an open subset $U\subset X$ containing the generic point of $W$ 
such that  $(K_U+\Delta|_U)$ is $\r$-Cartier (hence lc) and $W\cap U$ is an 
lc center of $(U, \Delta|_U)$.

Let $f:X\to S$ be a proper morphism and $L$ 
an $\r$-Cartier, $f$-nef divisor on $X$. Then $L$   is called {\it log $f$-big} if 
$L|_W$ is   big on the generic fiber of $f|_W:W\to f(W)$ for every lc center $W$ of $(X, \Delta)$ and also for every irreducible component $W\subset X$. Note that this notion depends on $\Delta$. 

If $0\leq \Delta'\leq\Delta$ and $(X, \Delta)$ is potentially lc then so is
$(X, \Delta')$ and every lc center of $(X, \Delta')$ is also an lc center of $(X, \Delta)$. Thus if $L$ is  log $f$-big on $(X, \Delta)$ then it is also log $f$-big on $(X, \Delta')$.

{\it Note.} There is another sensible way to define an ``lc center'' of a
potentially lc pair  $(X, \Delta)$ as an irreducible subvariety $W\subset X$ such that  $W\cap U$ is an
lc center of  $(U, \Delta|_{U}+\Theta_U)$ for every open subset $U\subset X$ and for every  effective divisor $\Theta_U$ such that $(U, \Delta|_{U}+\Theta_U)$ is lc (and $W\cap U$ is nonempty). With this definition, there are more ``lc centers.''

As a typical example, set  $Q:=(xy=uv)\subset \a^4$,
$A_1:=(x=u=0)$ and $A_2:=(y=v=0)$. Then
$\bigl(Q, A_1+A_2\bigr)$ is potentially lc. Its lc centers are the  divisors $A_1, A_2$. However, if $\Delta$ is any effective divisor such that  $\bigl(Q, A_1+A_2+\Delta\bigr)$ is lc then
the origin is also an lc center; cf.\  \cite[4.8]{ambro}, 
\cite[4.41]{kk-singbook} or \cite[6.3.11]{fuj-book}.

For our current purposes the first variant works better.  
\end{defn}

The following is proved in \cite{ambro} and \cite[1.10]{MR3238112}, see also 
\cite[Sec.5.7]{fuj-book} and \cite[6.3.5]{fuj-book}, where it is called a Reid-Fukuda--type vanishing theorem.

\begin{thm}[Ambro-Fujino vanishing theorem]
 \label{fujino.1.10.thm}   Let $(X, \Delta)$ be an slc pair and $D$ a Mumford $\z$-divisor on $X$.  Let $f:X\to S$ be a proper morphism. 
Assume that $D\simr K_X+L+\Delta$, where $L$ is $\r$-Cartier, $f$-nef and  log $f$-big.  Then
$$
R^if_*\o_X(D)=0\qtq{for} i>0. \qed
$$ 
\end{thm}

Combining it with Proposition~\ref{pot.lc.can.mod} gives the following variant, but only for normal varieties.

\begin{cor} \label{FKV.cor} Let $(X, \Delta)$ be a potentially lc pair, $f:X\to S$ a proper morphism  and $D$ a Weil $\z$-divisor on $X$.
Assume that $D\simr K_X+L+\Delta$, where $L$ is $\r$-Cartier, $f$-nef and  log $f$-big.  Then
$R^if_*\o_X(D)=0$ for $ i>0$. 
\end{cor}

Proof. Let $\pi: (X^c, \Delta^c)\to (X, \Delta)$ be as in Proposition~\ref{pot.lc.can.mod}. Then
$D^c\simr K_{X^c}+\pi^*L+\Delta^c$. Since  $\pi$ is a local isomorphism at every lc center of $(X^c, \Delta^c) $, we see that $\pi^*L$ is log $f\circ\pi$-big on
$(X^c, \Delta^c) $. By  Theorem~\ref{fujino.1.10.thm} this implies that
$$
\begin{array}{l}
R^i(f\circ\pi)_*\o_{X^c}(D^c)=0 \qtq{for} i>0\qtq{and}\\
R^i\pi_*\o_{X^c}(D^c)=0 \qtq{for} i>0.
\end{array}
$$
The Leray spectral sequence  now gives that $R^if_*\o_X(D)=0$ for $i>0$. \qed

\begin{say}[Log plurigenera]\label{log.pl.gen.say}
On an slc pair  $(X, \Delta)$, the best analogs of the multiples of the canonical divisor  are the divisors of the form
$mK_{X}+\rdown{m\Delta}$. While it is probably not crucial, it is convenient to know when their higher cohomologies vanish. Even if $K_X+\Delta$ is $\q$-Cartier and ample, the divisors $mK_{X}+\rdown{m\Delta}$ need not be
$\q$-Cartier and, even if $\q$-Cartier, need not be ample. 
In order to understand the situation, note that
$$
mK_{X}+\rdown{m\Delta}\simr
 K_{X}+(m-1)(K_{X}+\Delta) +\rdown{m\Delta}-(m-1)\Delta.
\eqno{(\ref{log.pl.gen.say}.1)}
$$
If $K_{X}+\Delta$ is nef and big then  (\ref{log.pl.gen.say}.1)  has the expected form for vanishing theorems, except that
$\rdown{m\Delta}-(m-1)\Delta $ need not be effective and 
the whole divisor need not be $\q$-Cartier.

If $\coeff\Delta\subset \{\frac12, \frac23, \frac34, \dots, 1\}$ then 
 $0\leq \rdown{m\Delta}-(m-1)\Delta\leq \Delta $ for every $m$ and 
 the second problem  is remedied by an application of Corollary~\ref{FKV.cor}. This leads to the following. 

\end{say}

\begin{prop} \label{stand.c.hi.coh.prop}
 Let $(X, \Delta)$ be an lc pair  with $\coeff\Delta\subset \{\frac12, \frac23, \frac34, \dots, 1\}$  and  $f:X\to S$ a proper morphism. Assume that
$K_X+\Delta$ is $f$-nef and log $f$-big.  Then
$$
R^if_*\o_{X}\bigl(mK_{X}+\rdown{m\Delta}\bigr)=0\qtq{for all} m\geq 2,\ i>0.
$$ 
\end{prop}

Proof.
If $m\geq 2$ then $(m-1)(K_{X}+\Delta)$ is $f$-nef and log $f$-big on
$(X,\Delta)$. Since $0\leq \rdown{m\Delta}-(m-1)\Delta\leq \Delta$, 
the pair $\bigl(X, \rdown{m\Delta}-(m-1)\Delta\bigr)$ is potentially lc and 
$(m-1)(K_{X}+\Delta)$ is also log $f$-big on
$\bigl(X,\rdown{m\Delta}-(m-1)\Delta\bigr)$. Thus
Corollary~\ref{FKV.cor} shows the required vanishing. \qed
\medskip

\begin{say}[Proof of Corollary \ref{st.pg.inv.cor}] By Theorem~\ref{half.pgsh.inv.thm} the sheaves $ \omega_{X/S}^{[m]}\bigl(\rdown{m\Delta}\bigr)$ are flat over $S$ and commute with base change. 
If all the fibers are normal then
the higher cohomologies of the fibers vanish by Proposition~\ref{stand.c.hi.coh.prop}. Thus  (\ref{st.pg.inv.cor}.1--2) are implied by the Cohomology and Base Change theorem. 

In general, pick $s\in S$,  let 
 $T$ be the spectrum of a DVR and $h:T\to S$ a morphism mapping the  generic point of $T$  to a generic point of $S$ and the closed point of $T$ to $s$. 
After pull-back we get $f_T: (X_T, \Delta_T)\to T$. Here  we can use Proposition~\ref{stand.c.hi.coh.prop} and conclude that 
$$
H^i\bigl(X_s, \o_{X_s}\bigl(mK_{X_s}+\rdown{m\Delta_s}\bigr)\bigr)=0\qtq{for all} m\geq 2,\ i>0.
$$
Thus the previous argument applies even if there are some non-normal fibers.
\qed
\end{say}


\section{Examples}\label{sec.4}

First we give some examples of  small modifications of lc singularities.

\begin{exmp} \label{sm.mod.CY.cnes.exmp}

(\ref{sm.mod.CY.cnes.exmp}.1) Let $(x\in X)$ be a cone over an Abelian variety $A$.
We see next that   small modifications of  $(x\in X)$ correspond to Abelian subvarieties
$A'\subset A$. So there are no small modifications if $A$ is simple but there are infinitely many small modifications if $A=E^m$ for some elliptic curve $E$ and $m\geq 2$.

Indeed, let $\pi:Y\to X$ be the cone-type resolution with exceptional divisor $F\cong A$. Note that $(Y, F)$ is canonical and $F$ is the unique divisor over $X$ with negative discrepancy.
 Let $\tau:Z\to X$ be a  small modification and  $Y'\to Z$  a $\q$-factorial  modification that extracts $F$; cf.\ \cite[1.38]{kk-singbook}. Then $Y'\map Y$ is an isomorphism in codimension 1, hence an isomorphism since $Y$ does not contain any proper rational curves, cf.\ \cite[VI.1.9]{rc-book}. Thus
we get a morphism with connected fibers $\tau_A:A\to \ex(\tau)$
and $Z$ is described by $\ker\tau_A$.

(\ref{sm.mod.CY.cnes.exmp}.2) Small modifications of cones over K3 surfaces can be obtained from  elliptic structures (contract the  cone-type resolution along the elliptic curves) and also from  configurations of $(-2)$-curves  (log-flop the curves and then contract the K3).  Thus there can be infinitely many different small modifications. The existence of  log flops between small modifications should follow from \cite{Kovacs94}. 

(\ref{sm.mod.CY.cnes.exmp}.3) An   example  showing that  log flops do not  always exist is given in \cite[96]{k-exrc}. Here is another one.

Let $P\subset \p^2$ be a set of 9 points in very general position. Set $S:=B_P\p^2$, let $E\subset S$ be the union of the 9 exceptional curves and   $C\subset S$  the birational transform of the unique cubic through $P$. Since the points are in very general position,
$H^0\bigl(S, \o_S(nC)\bigr)=1$ for every $n\geq 0$. 

Let $L\subset S$ be the birational transform of a general line; then  $C\sim 3L-E$ and $H:=C+L$ is ample. 

Let $X:=C_a(S, H)$ denote the affine cone over $S$ and   $D:=C_a(C, H|_C)\subset X$ the cone over $C$. Since  $K_S+C\sim 0$, the pair $(X, D)$ is lc.  (See \cite[3.1]{kk-singbook} for basic results on cones.)  Let $g:Y\to X$ be the cone-type resolution with exceptional divisor $F\cong S$. 

One can flop the 9 curves $E\subset F\subset Y$ and then contract 
the resulting $F'\cong \p^2$ to a  point (of type $\a^3/\tfrac14(1,1,1)$). This gives a small, lc modification  $\pi:X'\to X$ such that $D'$, the birational transform of $D$, is anti-ample. Thus $X'\cong \proj_X \tsum_{n=0}^{\infty} \o_X(-nD)$. 

We claim that $X'\to X$  does not have a log flop. The log flop  would be given by
$\proj_X \tsum_{n=0}^{\infty} \o_X(nD)$, but we check next that 
 this ring is not finitely generated. 

Since $X$ is a cone over $S$, 
$$
H^0\bigl(X,  \o_X(nD)\bigr)=\tsum_{m=0}^{\infty}H^0\bigl(S, \o_S(nC+mH)\bigr).
$$
Thus we need to show that 
$$
\tsum_{n,m=0}^{\infty}H^0\bigl(S, \o_S(nC+mH)\bigr)
$$ is not finitely generated.
The problem is with the $m=1$ summands. All of the $ \o_S(nC+H)$ are globally generated,
but  the image of
$$
 H^0\bigl(S, \o_S((n-r)C)\bigr)\times H^0\bigl(S, \o_S(rC+H)\bigr)\to
H^0\bigl(S, \o_S(nC+H)\bigr)
$$
consists of sections that vanish along $C$. 

Note also that if the 9 points $P$ are the base points of a cubic pencil, 
then $|C|$ is a base-point free elliptic pencil. The corresponding contraction $Y\to X''\to X$ gives the log flop of $X'\to X$.   



\end{exmp}

The following relates to Question~\ref{tech.prop.B.ques}.  

\begin{exmp}\label{2/3.v2.exmp} Set $Q:=(xy=uv)\subset \a^4$ and let
$|A|$ and $|B|$ denote the 2 families of planes on $Q$. Fix $n\geq 2$ and consider
the klt pair. 
$$
\bigl(Q, \Delta:=\bigl(\tfrac23-\tfrac{1}{3n-1}\bigr)(A_1+A_2+A_3)+
\bigl(1-\tfrac{1}{n}\bigr)B_1+\bigl(1-\tfrac{1}{n(3n-1)}\bigr)B_2\bigr).
$$
By direct computation we see that
$$
\omega_Q^{[3n]}\bigl(\rdown{3n\Delta}\bigr)\cong \o_Q(-6A-4B).
$$
This has only depth 2 at the origin  by \cite[3.15.2]{kk-singbook}.
\end{exmp}

Next we describe all CM divisorial sheaves on certain singularities and see how this compares with the conclusions of Proposition~\ref{tech.prop.0} and 
Theorem \ref{half.pgsh.inv.thm}. (Note that for a 3-dimensional isolated singularity $(x, X)$ a divisorial sheaf $F$ is CM  $\Leftrightarrow$ $F$ is $S_3$
 $\Leftrightarrow$ $\depth_xF\geq 3$.)

\begin{exmp} \label{half.main.exmp}
Start with $Q:=(xy-uv=0)\subset \a^4$ and 
the divisors   $B_0:=(y=v=0)$, $B:=(x=u=0)$ and $T:=(v=u^n)$. 
Then $T\cong (xy-u^{n+1}=0)$.

Next take  quotient by a $\mu_n$-action
$X:=Q/\frac1{n}(1,0,1,0)$. The quotient map
$\tau:Q\to X$ ramifies along $B$.
In $X$ consider the (reduced)  divisors  $D_0:=\tau(B_0)$, $D:=\tau(B)$ and $S:=\tau(T)$. Note that  $(x/v)^n$ is a rational function on $X$ and
it equals $(u/y)^n$. This shows that $D\sim nD_0$ and
we compute that
$$
\cl(X)=\z[D_0] \qtq{and} K_X\sim -(n-1)D_0.
\eqno{(\ref{half.main.exmp}.1)}
$$
Furthermore, 
$\tau^*(K_X+S+aD_0+bD)=K_Q+T-(n-1)B+aB_0+nbB$  is $\q$-Cartier iff $a+nb-n+1=0$. 
Using  \cite{r-c3f} or \cite[2.43]{kk-singbook} and the convexity of lc boundaries gives the following.
\medskip

{\it Claim \ref{half.main.exmp}.2.} For any $0\leq \lambda\leq 1$ the pair
$$
\bigl(X, S+\lambda\bigl(D_0+\tfrac{n-2}{n}D\bigr)+(1-\lambda)\tfrac{n-1}{n}D\bigr)
\qtq{is lc.} \qed
$$
We can also get a complete description of all CM divisorial sheaves.

\medskip

{\it Claim \ref{half.main.exmp}.3.}  $\o_X(mD_0)$ is  CM iff $-n\leq m\leq 1$.
\medskip

Proof. $Q$ is a cyclic cover of $X$ ramified along $D\sim nD_0$, thus
$$
\pi_*\o_Q\cong \o_X+\o_X(-D_0)+\cdots+ \o_X(-(n-1)D_0).
$$
So $\o_X(mD_0)$ is  CM
for $-(n-1)\leq m\leq 0$.
For $m=-n$ we use that  $\o_X(-nD_0)\cong \o_X(-D)$
which is the $\mu_n$-invariant part of $\o_Q(-B)$ hence CM, cf.\ \cite[3.15.2]{kk-singbook}. 
By Serre duality  we get that $\o_X(D_0)$ is CM.

For $m<-n$ set $r:=-m-n$. Note that $D+rD_0$ is not $S_2$, since the cokernel of
$\o_{D+rD_0}\to \o_{D}+\o_{rD_0}$  is supported at the origin.
Then the exact sequence
$$
0\to \o_X(-D-rD_0)\to \o_X\to \o_{D+rD_0}\to 0
$$
shows that $\o_X(mD_0)$ is not CM, cf.\ \cite[2.60]{kk-singbook}. 
 By Serre duality this gives that $\o_X(mD_0)$ is not CM for $m\geq 2$. \qed
\medskip

{\it Corollary \ref{half.main.exmp}.4.} If $n\geq 3$ then
$\omega_X^{[m]}\cong \o_X\bigl(-m(n-1)D_0\bigr)$ is CM only for $m=0,1$. \qed
\medskip

Next we see how this example compares with the conclusion of
Theorem \ref{half.pgsh.inv.thm} for various choices of the boundary $\Delta$.
\medskip

{\it Example \ref{half.main.exmp}.5.} Assume that $n$ is even and set  $\Delta:=\frac{n-1}{n+1}(D_0+D)$. Then $K_X+\Delta\simq 0$ and
$(X, S+\Delta)$ is lc by  (\ref{half.main.exmp}.2). 

Observe that 
 $\frac{n+2}{2}K_X+\rdown{\frac{n+2}{2}\Delta}\simq -nD_0$, thus  Proposition~\ref{tech.prop.B} implies that $\o_X(-nD_0)$ is  CM. This coincides with the  lower bound in 
(\ref{half.main.exmp}.3). 
We also see that 
 $$
mK_X+\rdown{r\Delta}\simq \bigl(-m(n-1)+(n+1)\rdown{\tfrac{r(n-1)}{n+1}}\bigr)D_0.
$$
Elementary estimates  show that    $\omega_X^{[m]}\bigl(\rdown{r\Delta}\bigr)$ is not CM whenever $|m-r|\geq 2$ and $n\geq 4$. 
\medskip

{\it Example \ref{half.main.exmp}.6.}
Set  $\Delta:=\frac{n-1}{2n+1}(D_0+D'+D'')$ where $D',D''\in |D, nD_0|$ are two general divisors. Then $K_X+\Delta\simq 0$ and
$(X, S+\Delta)$ is lc by (\ref{half.main.exmp}.2). The coefficients $\frac{n-1}{2n+1}$ converge to $\frac12$ from below. 

Take $m=2$, Then $2K_X+\rdown{2\Delta}\simq -(2n-2)D_0$, thus
$\omega^{[2]}_X\bigl(\rdown{2\Delta}\bigr)=\omega^{[2]}_X$ is  not CM for $n\geq 3$ by (\ref{half.main.exmp}.3).  This shows that the bound $\frac12$ is sharp in  Proposition~\ref{tech.prop.B} and in Theorem~\ref{half.pgsh.inv.thm}.

\medskip

{\it Example \ref{half.main.exmp}.7.} First set $\Delta:=\frac{n-1}{n+1}(D_0+D)$,  $B:=-D_0-D$ and $\Delta_3:=D_0+D$. Then
$\Delta_3\leq \rup{\Delta}$ yet $\o_X(B)$ is not CM. Thus the condition $\rdown{\Delta_3}\leq \rdown{\Delta}$ is necessary in  Proposition~\ref{tech.prop.2}.

Next take   $\Delta:=\frac12(D'_0+D''_0)+\frac{n-2}{n}D$ where $D'_0,D''_0\in |D_0|$ are two general divisors.  Then $K_X+\Delta\simq 0$ and
$(X, S+\Delta)$ is lc by (\ref{half.main.exmp}.2). 
Take 
$$
B:=-D_0-D \qtq{and}
\Delta_3:= \tfrac{1+\epsilon}{2}(D'_0+D''_0)+\tfrac{n-\epsilon}{n}D.
$$ Then  $\rdown{\Delta_3}=0$, 
$\Delta_3\leq\rup{\Delta}$ but $\o_X(B)$ is not CM.   This shows that 
$\Delta^{=1/2}$ needs special  handling in  Proposition~\ref{tech.prop.2}.

\medskip

{\it Example \ref{half.main.exmp}.8.} Let $0<\epsilon\ll 1$ be irrational and  set  $$\Delta:=\tfrac{n-1}{n+1}\bigl((1+n\epsilon)D_0+(1-\epsilon)D\bigr).$$ Then $K_X+\Delta\simr 0$ and
$(X, S+\Delta)$ is lc by  (\ref{half.main.exmp}.2). Also,
$$
mK_X+\rdown{m\Delta}\simr -\Bigl(\bigl\{m\tfrac{n-1}{n+1}(1+n\epsilon)\bigr\}+
n\bigl\{m\tfrac{n-1}{n+1}(1-\epsilon)\bigr\}\Bigr)D_0,
$$
and the right hand side is always a strictly negative multiple of $D_0$ for $m\neq 0$.
Thus $mK_X+\rdown{m\Delta}$ is not Cartier for every  $m\neq 0$.

\end{exmp}

The next example shows that in Proposition~\ref{stand.c.hi.coh.prop} it is not enough to assume that the coefficients of $\Delta$ are close to 1. 
In fact, vanishing fails for every other value of the coefficients.

\begin{exmp} \label{c.near.1.h.coh.exmp} Set $X=\p^n$ and choose  $0<\epsilon<\frac1{n+2}$. Let $H_i$ be hyperplanes in general position and set 
$$
\Delta=
\bigl(1-\tfrac{1+\epsilon}{n+3}\bigl)\bigl(H_1+\cdots+H_{n+2}\bigr).
$$
Then $K_X+\Delta\simr \tfrac{1-(n+2)\epsilon}{n+3} H$ is ample and
$$
\rdown{(n+3)\Delta}=(n+1)\bigl(H_1+\cdots+H_{n+2}\bigr).
$$
Therefore
$$
(n+3)K_{\p^n}+\rdown{(n+3)\Delta}\simr K_{\p^n},
$$
which has nonzero higher cohomology.

Finally note that if $\epsilon'=\frac1{n+2} $ then  $1-\tfrac{1+\epsilon'}{n+3}=1-\tfrac{1}{n+2}$. Thus as $\epsilon$ varies in the interval  $\bigl(0, \frac1{n+2}\bigr)$,  the  values of
$1-\tfrac{1+\epsilon}{n+3}$ cover the open interval
$\bigl(1-\tfrac{1}{n+2}, 1-\tfrac{1}{n+3}\bigr)$. 
\end{exmp}

Each of the above examples gives many more in one dimension higher.

\begin{exmp} \label{c.near.1.h.coh.exmp.2}
Let  $Y'$ be a smooth projective variety of dimension $n+1$. Pick a point 
$y'\in Y'$ and let  $H'_1, \cdots, H'_{n+2}$ be general smooth divisor passing through $y'$. 
Let $\pi:Y\to Y'$ be the blow-up of $y'$ with exceptional divisor $E$. Let $H_i$ denote the birational transform of $H'_i$.  Our example is
$$
(Y, E+\Delta)\qtq{where}
\Delta=
\bigl(1-\tfrac{1+\epsilon}{n+3}\bigl)\bigl(H_1+\cdots+H_{n+2}\bigr).
$$
Note that $(E, \diff_E\Delta)$  is isomorphic to the example in
(\ref{c.near.1.h.coh.exmp}). In partcular, 
our previous computations show that
$$
R^n\pi_* \omega_Y^{[n+3]}\bigl((n+3)E+ \rdown{(n+3)\Delta}\bigr)\cong k(y').
$$
If the $H_i$ are chosen sufficienty ample then we get that
$$
h^n\bigl(Y, \omega_Y^{[n+3]}\bigl((n+3)E+ \rdown{(n+3)\Delta}\bigr)\bigr)=1.
$$
Note that such examples appear very naturally if we try to compactify the moduli space of $n+1$ dimensional varieties with $n+2$ divisors on them. 
\end{exmp}

\begin{exmp} \label{c.near.1.h.coh.exmp.3} Fix integers $n\geq m\geq 2$. Set ${\mathbf P}=\p^{mn-n-1}$ and choose  $0<\epsilon<\frac1{mn}$. Let $H_0, \dots, H_{mn}$ be hyperplanes in general position and set 
$$
\Delta=H_0+
\bigl(1-\tfrac1{m}-\epsilon\bigl)\bigl(H_1+\cdots+H_{mn}\bigr).
$$
Then $K_{\mathbf P}+\Delta\simr (1-mn\epsilon) H$ is ample and
$$
\rdown{m\Delta}=mH_0+(m-2)\bigl(H_1+\cdots+H_{mn}\bigr).
$$
Therefore
$$
mK_{\mathbf P}+\rdown{m\Delta}\simr K_{\mathbf P}-(n-m)H,
$$
which has nonzero  cohomology in degree $mn-n-1$ if $n\geq m$. 
\end{exmp}

\begin{exmp} \label{c.near.1.h.coh.exmp.4}
Fix integers $n\geq m\geq 2$ and  choose  $0<\epsilon<\frac1{mn}$.
Let  $Y'$ be a smooth projective variety of dimension $mn-n-1$. Pick a point 
$y'\in Y'$ and let  $H'_1, \cdots, H'_{mn}$ be general smooth divisor passing through $(y',0)\in Y'\times \a^1$. 
Let $\pi:X\to Y'\times \a^1$ be the blow-up of $(y',0)$ with exceptional divisor $E$. Let $H_i$ denote the birational transform of $H'_i$.  
Set $$
\Delta=
\bigl(1-\tfrac1{m}-\epsilon\bigl)\bigl(H_1+\cdots+H_{mn}\bigr)$$
and let $f:(X,\Delta)\to \a^1$ be the composite of the coordinate projection with $\pi$. 

If the $H_i$ are sufficiently ample then  $f:(X,\Delta)\to \a^1$ is a stable morphism and  $(X_t, \Delta_t)$ is a simple normal crossing pair for $t\neq 0$.  The fiber over the origin has 2 irreducible components.
One is $E$ and the other is $Y_0$, which  is the blow-up of $y'\in Y'_0$. 

Note that $\bigl(E, \diff_E(Y_0+\Delta)\bigr)$  is isomorphic to the example in
(\ref{c.near.1.h.coh.exmp.3}). In partcular, 
our previous computations show that
$$
R^{mn-n-1}\pi_* \omega_Y^{[m]}\bigl(\rdown{m\Delta}\bigr)\neq 0
$$
and the other higher direct images are 0. Thus,
if the $H_j$ are chosen sufficienty ample then we get that
$$
R^{mn-n-1}f_* \omega_Y^{[m]}\bigl(\rdown{m\Delta}\bigr)
$$
is a torsion sheaf supported at the origin $0\in \a^1$ 
and the other higher direct images are 0.
In this case the functions
$$
t\mapsto h^i\bigl(X_t, \omega_{X_t}^{[m]}(\rdown{m\Delta_t})\bigr)
$$
jump at $t=0$ for $i=mn-n-1$ and $i=mn-n-2$  but are constant for
other values of $i$. 
\end{exmp}


\def\cprime{$'$} \def\cprime{$'$} \def\cprime{$'$} \def\cprime{$'$}
  \def\cprime{$'$} \def\cprime{$'$} \def\dbar{\leavevmode\hbox to
  0pt{\hskip.2ex \accent"16\hss}d} \def\cprime{$'$} \def\cprime{$'$}
  \def\polhk#1{\setbox0=\hbox{#1}{\ooalign{\hidewidth
  \lower1.5ex\hbox{`}\hidewidth\crcr\unhbox0}}} \def\cprime{$'$}
  \def\cprime{$'$} \def\cprime{$'$} \def\cprime{$'$}
  \def\polhk#1{\setbox0=\hbox{#1}{\ooalign{\hidewidth
  \lower1.5ex\hbox{`}\hidewidth\crcr\unhbox0}}} \def\cdprime{$''$}
  \def\cprime{$'$} \def\cprime{$'$} \def\cprime{$'$} \def\cprime{$'$}
\providecommand{\bysame}{\leavevmode\hbox to3em{\hrulefill}\thinspace}
\providecommand{\MR}{\relax\ifhmode\unskip\space\fi MR }
\providecommand{\MRhref}[2]{%
  \href{http://www.ams.org/mathscinet-getitem?mr=#1}{#2}
}
\providecommand{\href}[2]{#2}

\bigskip

\noindent  Princeton University, Princeton NJ 08544-1000

{\begin{verbatim} kollar@math.princeton.edu\end{verbatim}}

\end{document}